\documentclass[11pt,a4paper,oneside]{amsart}
\usepackage[utf8]{inputenc}
\usepackage{mystyle}
\usepackage{csquotes}

\makeatletter
\@namedef{subjclassname@2020}{\textup{2020} Mathematics Subject Classification}
\makeatother

\newcommand{\termA}{\ensuremath{\ast}}
\title{Presheaves of groupoids as models for homotopy types}
\author{Léonard Guetta}
\address{Max~Planck~Institute~for~Mathematics,~Bonn,~Germany}
\email{guetta@mpim-bonn.mpg.de}
\urladdr{\href{http://guests.mpim-bonn.mpg.de/guetta/}{http://guests.mpim-bonn.mpg.de/guetta/}}
\date{\today}
\begin{document}
\maketitle
\begin{abstract}
  We introduce the notion of groupoidal (weak) test category, which is a small
  category $\sA$ such that the $\Grpd{}$\nbd{}valued presheaves over $\sA$
  models homotopy types in a ``canonical and nice'' way. The definition does not
  require \emph{a priori} that $\sA$ is a (weak) test category, but we prove two
  important comparison results: (1) every weak test category is a groupoidal
  weak test category, (2) a category is a test category \emph{if and only if} it
  is a groupoidal test category.

  As an application, we obtain new models for homotopy types, such as the
  category of groupoids internal to cubical sets with or \emph{without}
  connections, the category of groupoids internal to cellular sets, the category
  of groupoids internal to semi-simplicial sets, etc.

  We also prove, as a by-product result, that the category of groupoids internal
  to the category of small categories models homotopy types.
\end{abstract}
\tableofcontents

\section*{Introduction}
In his famous manuscript \emph{Pursuing Stacks} from 1983
\cite{grothendieck1983pursuing}, Grothendieck introduced the theory of
\emph{test categories}. Informally speaking, a test category is a small category
$\sA$ such that the category $\psh{\sA}$ of $\Set$\nbd{}valued presheaves over
$\sA$ models homotopy types in a canonical way, which was axiomatized by
Grothendieck. ``Models homotopy types'' means here that there is a particular
class of morphisms on $\psh{\sA}$, such that the localization of $\psh{\sA}$
with respect to this class of morphisms is equivalent to the category $\Hot$ of
\emph{homotopy types}, that is the category of CW\nobreakdash-complexes and
homotopy classes of continuous maps between them.\footnote{In practise, this
  equivalence of homotopy categories will always come from an equivalence at the
  level of $(\infty,1)$\nbd{}categories.} The archetypal example of a test
category is of course the category $\Delta$ of finite non-empty ordinals, for
which it has been known since the famous result of Milnor
\cite{milnor1957geometric} that the category $\psh{\Delta}$ of simplicial sets
models homotopy types.

Examples of test categories abound, such as the cubical category
\cite[Corollaire 8.4.13]{cisinski2006prefaisceaux}, the cubical category with
connections \cite{maltsiniotis2009categorie}, Joyal's $\Theta$ category
\cite{cisinski2011categorie}, the dendroidal category $\Omega$
\cite{ara2019dendroidal}, etc. The point of view of the theory of test
categories is that the category of presheaves on \emph{any} test category should
be, in some sense, as good a model of homotopy types as the category of
simplicial sets. But modeling homotopy types is not the only good homotopical
property of the category of simplicial sets. For example, we have the following
results:
\begin{rome}
\item the category $\Ab(\psh{\Delta})$ of abelian groups internal to simplicial
  sets models \emph{homology types}, that is chain complexes in non-negative
  degree up to quasi-isomorphisms\footnote{The classical ``Dold--Kan''
    equivalence even says that there is an equivalence of categories between
    $\Ab(\psh{\Delta})$ and the category of chain complexes in non-negative
    degree. However, from a homotopical point view, it is the homotopical result
    stated above which is relevant.} \cite{dold1958homology},
  \cite{kan1958functors},
\item the category $\Grp(\psh{\Delta})$ of groups internal to simplicial sets
  models pointed connected homotopy types \cite{kan1958homotopy},
\item the category $\Grpd(\psh{\Delta})$ of groupoids internal to simplicial
  sets models homotopy types \cite[Theorem 8.3]{crans1995quillen}, \cite[Theorem
  10]{joyal1996homotopy}.
\end{rome}
In this article, we focus on property (iii) above and prove its generalization,
whose precise formulation requires the basic setup of the theory of test
categories, which we recall now.
      
For any small category $\sA$, Grothendieck considers the following canonical
functor
\[
  i_{\sA} \colon \psh{\sA} \to \Cat,
\]
which sends an object $X$ of $\psh{\sA}$ to its category of elements $\sA/X$. By
a result attributed to Quillen by Illusie \cite[Corollaire
3.3.1]{illusie1972complexe}, the category $\Cat$ models homotopy types, when
equipped with the class $\cW_{\infty}$ of functors $u \colon \sC \to \sD$ that
induce a weak homotopy equivalence between the classifying spaces of $\sC$ and
$\sD$.
Using this, we can the define the \emph{homotopy type} of a presheaf $X$ over
$\sA$ as the category $i_{\sA}(X)=\sA/X$ seen as an object of
$\Cat[\cW_{\infty}^{-1}] \simeq \Hot$. This definition has a very natural
interpretation: the category $\sA/X$ is nothing but the so--called
\emph{Grothendieck construction} of $X$, and by a theorem of Thomason
\cite{thomason1979homotopy}, this is the homotopy colimit of $X$, seen as a
(discrete) space\nbd{}valued presheaf. More generally, we define a class of weak
equivalences in $\psh{\sA}$ as
\[
  \cW_{\psh{\sA}}:=i_{\sA}^{-1}(\cW_{\infty}),
\]
and consider the induced functor at the level of homotopy categories
\[
  \overline{i_{\sA}} \colon \psh{\sA}[\cW_{\psh{\sA}}^{-1}] \to
  \Cat[\cW_{\infty}^{-1}] \simeq \Hot.
\]

The category $\sA$ is a \emph{pseudo-test category} if this functor is an
equivalence of categories, a \emph{weak test category} if the right adjoint of
$i_{\sA}$ (which always exists) also preserves weak equivalences and induces an
equivalence of homotopy categories, and, finally, a \emph{test category} if it
is a weak test category and for every object $a$ of $\sA$, $\sA/a$ is a weak
test category.
       
What about $\Grpd(\psh{\sA})$ now? There are two trivial but essential
observations: (1) the category of groupoids internal to $\psh{\sA}$ is
equivalent to the category $[\sA^{\op},\Grpd]$ of presheaves over $\sA$ with
values in the category of groupoids, (2) the Grothendieck construction is also
defined for $\Grpd$\nbd{}valued presheaves, and so we can define a functor
\[
  I_{\sA} \colon \Grpd(\psh{\sA}) \to \Cat,
\]
where $I_{\sA}(X)$ is the Grothendieck construction of $X$. The generalization
of the situation for $\psh{\sA}$ described earlier is then straightforward. We
define a canonical class of weak equivalences on $\Grpd(\psh{\sA})$ as
\[
  \cW_{\Grpd(\psh{\sA})}:=I_{\sA}^{-1}(\cW_{\infty}),
\]
and consider the functor induced at the level of homotopy categories
\[
  \overline{I_{\sA}} \colon \Grpd(\psh{\sA})[\cW_{\Grpd(\psh{\sA})}^{-1}] \to
  \Cat[\cW_{\infty}^{-1}] \simeq \Hot.
\]
Then, we can define the notions of \emph{groupoidal pseudo-test category},
\emph{groupoidal weak test category} and \emph{groupoidal test category} as
perfect analogues of the usual notions for $\Set$\nbd{}valued
presheaves.\footnote{We will also define the notions of ``groupoidal local test
  category'' and ``groupoidal strict test category''.} In fact, almost all of
the results from the usual theory also work in the groupoidal theory, once they
have been appropriately adapted. Once the new theory is well established, we can
compare the $\Set$\nbd{}valued notions with the $\Grpd$\nbd{}valued notions and
we obtain the main result of this paper.
\begin{theoremintro}\label{thm:intro}
  Let $\sA$ be a small category.
  \begin{enumerate}[label=(\roman*)]
  \item $\sA$ is a groupoidal test category \emph{if and only if} it is a test
    category,
  \item if $\sA$ is a weak test category, then it is also groupoidal weak test
    category.
  \end{enumerate}
\end{theoremintro}
Before getting to the immediate applications of this theorem, let us comment on
an important point.

Part $(ii)$ of the previous theorem, which implies that for a (weak) test
category $\sA$, the category $\Grpd(\psh{\sA})$ models homotopy types, does not
really come as a surprise, and here is a sketch of a very short proof of this
mere fact (see \cref{section:grpdcatmodel} for details). One can easily show
that if we equip the category $\Grpd(\Cat)$ of groupoids internal to the
category of small categories with the class of weak equivalences induced by the
diagonal of a well-defined bisimplicial nerve, then for \emph{every} weak test
category $\sA$, we have an equivalence of homotopy categories\footnote{Note
  however that this does \emph{not} work if we replace (internal) groupoids by
  (internal) groups. Indeed, simplicial groups models pointed connected homotopy
  types, but groups internal to categories (i.e.\ crossed modules) only model
  pointed connected homotopy $2$\nbd{}types (cf.\ \cref{rem:grpdcatmodel}).}
\[
  \Ho(\Grpd(\psh{\sA}))\simeq\Ho(\Grpd(\Cat)).
\]
Then, by applying this isomorphism twice, if $\sB$ is any other weak test
category, we have
\[
  \Ho(\Grpd(\psh{\sA}))\simeq \Ho(\Grpd(\psh{\sB})).
\]
In particular, for $\sB=\Delta$, we already know that $\Grpd(\psh{\Delta})$
models homotopy types, hence the desired result.

Nevertheless, the approach taken in this paper is \emph{completely} different
and the proof of \cref{thm:intro} comes as an easy consequence of the theory of
groupoidal test categories thoroughly developed beforehand. The advantage is at
least twofold:
\begin{enumerate}
\item we do not use the previously known case of $\Delta$, and hence we obtain a
  new (and much simpler) proof that $\Grpd(\psh{\Delta})$ models homotopy types,
\item we also obtain part $(i)$ of \cref{thm:intro}, which gives a converse for
  groupoidal test categories: if $\sA$ is a groupoidal test category,
  \emph{then} $\sA$ is a test category.
\end{enumerate}
This second point is rather astonishing; it says that the fact that
$\Grpd(\psh{\sA})$ models homotopy types (in a canonical way) is a necessary
\emph{and} sufficient condition so that $\psh{\sA}$ models homotopy types. For
the sake of comparison, let us consider the category $\Cat(\psh{\sA})$ of
category objects in $\psh{\sA}$. It is an easy exercise to show that this
category models homotopy types (in a canonical way) if and only the nerve of
$\sA$ is weakly contractible. In particular, it is \emph{not} a sufficient
condition to ensure that $\psh{\sA}$ models homotopy types. Somehow, the
homotopy theory of $\Grpd{}$\nbd{}valued presheaves is unexpectedly well
behaved.

As an application of \cref{thm:intro}, 
we obtain a plethora of new models of homotopy types:
\begin{itemize}[label=-]
\item the category $\Grpd(\psh{\square})$ of groupoids internal to cubical sets
  with or without connections,
\item the category $\Grpd(\psh{\Theta})$ of groupoids internal to cellular sets,
\item the category $\Grpd(\psh{\Omega})$ of groupoids internal to dendroidal
  sets,
\item the category $\Grpd(\psh{\Delta'})$ of groupoids internal to
  semi-simplicial sets,
\item etc.
\end{itemize}
The example of cubical sets \emph{without} connections is certainly worth
noticing. Indeed, it is a folkloric result that cubical groups without
connections are not Kan in general (see~\cite{tonks1992cubical} for example),
and, consequently, it is believed that the category of cubical groups without
connections does \emph{not} model pointed connected homotopy types, contrary to
simplicial groups. It can then come as a surprise that the category of groupoids
internal to cubical sets without connections models homotopy types. Once again,
$\Grpd$\nbd{}valued presheaves seem to behave particularly well and hopefully
this justifies their study.
      
Let us end this introduction with a quick word on what is \emph{not} treated in
the present article. In his book \cite{cisinski2006prefaisceaux}, Cisinski
showed that the category of ($\Set$\nbd{}valued) presheaves on any test category
admits a model structure where the weak equivalences are the ones canonically
defined by Grothendieck. The generalization of this for $\Grpd$\nbd{}valued
presheaves (known when $\sA=\Delta$~\cite{crans1995quillen,joyal1996homotopy})
is not addressed at all here and is left as future work.

\begin{named}[Organization of the paper]The first section is a preliminary
  section recalling some basic homotopical algebra needed in the rest of the
  paper. The second section is a quick recollection, without proofs, of the
  classical theory of test categories. It is in the third section that we
  finally dive into the subject and give the basic setup of the homotopy theory
  of $\Grpd$\nbd{}valued presheaves, mimicking Grothendieck's axiomatic for the
  homotopy theory of $\Set$\nbd{}valued presheaves. In the fourth, fifth and
  sixth sections, we respectively develop the theory of groupoidal weak test
  categories, groupoidal test categories and groupoidal strict test categories.
  The seventh section is dedicated to the comparison of the classical theory
  with the groupoidal theory and in particular we obtain the main theorem stated
  previously in the introduction. Then, in the eighth section, we give an
  alternative definition of the weak equivalences of $\Grpd$\nbd{}valued
  presheaves, which makes the link with previous existing works on simplicial
  groupoids. Finally, the ninth section is a ``bonus'' and almost independent
  section from the rest of the paper, where we prove that $\Grpd(\Cat)$ model
  homotopy types.
\end{named}
\begin{named}[Acknowledgements]
  I am very thankful to Georges Maltsiniotis for the countless helpful
  conversations. The original idea for this paper stemmed from one of them. I am
  also grateful to the Max Planck Institute for Mathematics in Bonn for its
  hospitality and financial support while this article was being written.
\end{named}
      
\section{Preliminaries}
\begin{paragr}
  A \emph{category with weak equivalences} is a pair $(\cC,\cW)$, where $\cC$ is
  a category and $\cW$ is a class of morphisms of $\cC$, generically referred to
  as the \emph{weak equivalences}, which contains all isomorphisms and satisfy
  the 2-out-of-3 property. We say that $\cW$ is \emph{weakly saturated} if in
  addition it satisfies the following closure property: If $i \colon X \to Y$
  and $r \colon Y \to X$ are morphisms of $\cW$ such that $ri=\id_Y$ and $ir \in
  \cW$, then $r \in \cW$ (and thus so is $i$ by 2-out-of-3).

  When $\cC$ has pullbacks, we say that a morphism $f \colon X \to Y$ is a
  \emph{universal weak equivalence}, if for every pullback square
  \[
    \begin{tikzcd}
      X' \ar[r] \ar[d,"f'"'] & X \ar[d,"f"]\\
      Y' \ar[r] & Y, \ar[from=1-1,to=2-2,"\lrcorner",phantom,very near start]
    \end{tikzcd}
  \]
  the morphism $f'$ is a weak equivalence (in particular $f$ is a weak
  equivalence).

  By \emph{homotopy category} of a category with weak equivalences $(\cC,\cW)$,
  we mean the localization of $\cC$ with respect to $\cW$
  \cite{gabriel1967calculus}. It is denoted by $\Ho_{\cW}(\cC)$, or simply
  $\Ho(\cC)$ when there is no risk of confusion.
\end{paragr}
\begin{paragr}\label{def:homotopicalequiv}
  Let $(\cC,\cW)$ and $(\cC',\cW')$ be two categories with weak equivalences. A
  functor $F \colon \cC \to \cC'$ is said to \emph{preserve weak equivalences}
  if $F(\cW) \subseteq \cW'$. In this case, it induces a functor at the level of
  homotopy categories
  \[
    \overline{F} \colon \Ho(\cC) \to \Ho(\cC').
  \]
  A \emph{homotopy inverse} of such a functor, is a functor $G \colon \cC' \to
  \cC$ which preserves weak equivalences and such that there exists zigzags of
  natural transformations $FG \leftrightsquigarrow \id_{\cC'}$, and $GF
  \leftrightsquigarrow \id_{\cC}$, which are pointwise weak equivalences (see
  (ii) below). In this case, $F$ induces an equivalence of categories
  $\Ho(\cC)\simeq \Ho(\cC')$ and similarly for $G$.\footnote{In fact, $F$ and
    $G$ are then even Dwyer--Kan equivalences
    \cite{barwick2012characterization}, hence they induce an equivalence of
    $(\infty,1)$\nbd{}categories.}

  Finally, we say that an adjunction $\begin{tikzcd} L \colon \cC \ar[r,shift
    left] & \cC' \colon R \ar[l,shift left]\end{tikzcd}$ is a \emph{homotopical
    equivalence} if:
  \begin{rome}
  \item $L$ and $R$ preserve weak equivalences,
    \[
      L(\cW)\subseteq \cW' \text{ and } R(\cW') \subseteq \cW
    \]
  \item the unit and co--unit of the adjunction are pointwise weak equivalences,
    i.e.\ for every object $X$ of $\cC'$ and every object $Y$ of $\cC$ we have
    \[
      \varepsilon_X \colon LR(X) \to X \in \cW' \text{ and } \eta_Y \colon Y \to
      RL(Y) \in \cW.
    \]
  \end{rome}
  Note that in this case $L$ and $R$ are homotopy inverses to each other.
\end{paragr}
The following lemma is a very useful criterion to detect homotopical
equivalences.
\begin{lemma}\label{lemma:crithomotopicalequiv}
  Let $(\cC,\cW)$ and $(\cC',\cW')$ be two categories with weak equivalences
  and 
  $\begin{tikzcd} L \colon \cC \ar[r,shift left] & \cC' \colon R \ar[l,shift
    left]\end{tikzcd}$ an adjunction. If $\cW=L^{-1}(\cW')$, then $L \dashv R$
  is a homotopical equivalence if and only if the co--unit of the adjunction is
  a pointwise equivalence.

\end{lemma}
\begin{rem}
  The dual of the previous lemma is also true, but we won't need it in this
  paper.
\end{rem}
\begin{proof}
  First, let's prove that $R$ preserves weak equivalences. Let $f \colon X \to
  Y$ be a morphism a $\cC'$ that belongs to $\cW'$. Since $\cW=L^{-1}(\cW')$, we
  need to show that $LR(f) \in \cW'$ and this follows from the 2-out-of-3
  property of $\cW'$, the commutativity of the square
  \[
    \begin{tikzcd}
      LR(X) \ar[r,"LR(f)"] \ar[d,"\varepsilon_X"']& LR(Y) \ar[d,"\varepsilon_Y"] \\
      X \ar[r,"f"'] & Y,
    \end{tikzcd}
  \]
  and the fact that the co--unit is a pointwise weak equivalence.

  Now let's prove that the unit of the adjunction is also a pointwise weak
  equivalence. Let $Y$ be an object of $\cC'$. Since $\cW=L^{-1}(\cW')$, we need
  to show that $L(\eta_Y)$ is a weak equivalence. By the triangle identity, the
  following triangle is commutative
  \[
    \begin{tikzcd}
      L(Y) \ar[r,"L(\eta_Y)"] \ar[rr,"\id_{L(Y)}"',bend right] & LRL(Y)
      \ar[r,"\varepsilon_{L(Y)}"] & L(Y),
    \end{tikzcd}
  \]
  hence the desired result follows from the 2-out-of-3 property of $\cW'$ and
  the fact the co--unit is a pointwise weak equivalence.
\end{proof}
\section{Very brief recollection of the theory of test categories}
\emph{The goal of this section is only to provide a quick summary of the basic
  notions and results (without proofs) of the theory of test categories. For a
  detailed exposition, we refer the reader to Maltsiniotis' book on the subject
  \cite{maltsiniotis2005theorie}. }
\begin{notation}
  For a small category $\sA$, we denote by $\psh{\sA}$ the category of
  $\Set$\nbd{}valued presheaves over $\sA$, that is, the category of functors
  $\sA^{\op} \to \Set$ and natural transformations between them.

  We denote by $\Cat$ the category of small categories and functors between
  them. We use the notation $e$ for the terminal object of $\Cat$, that is, the
  category with one object and no non-identity morphism.

  For a small category $\sA$ and an object $a$ of $\sA$, we denote by $\sA/a$
  the slice category of $\sA$ over $a$. Explicitly, $\sA/a$ is the category
  whose objects are pairs $(a', p \colon a' \to a)$, where $a'$ and $p$ are
  respectively an object and a morphism of $\sA$, and whose morphisms $(a',p')
  \to (a'',p'')$ are morphisms $f \colon a' \to a''$ of $\sA$, such that
  $p''\circ f = p'$. Note that we have an obvious projection functor $\sA/a \to
  \sA$.

  More generally, if $u \colon \sA \to \sB$ is a morphism of $\Cat$ and $b$ is
  an object of $\sB$, we denote by $\sA/b$ the category whose objects are pairs
  $(a,q \colon u(a) \to b)$, where $a$ is an object of $\sA$ and $q$ is a
  morphism of $\sB$, and whose morphisms $(a,q) \to (a',q')$ are the morphisms
  $f \colon a \to a'$ of $\sA$ such that $q' \circ u(f) =q$. Note that we make
  the abuse of notation of not making $u$ appear in the notation $\sA/b$, but
  this category obviously depends on $u$.
\end{notation}
\begin{paragr}
  Let $\Delta$ be the category whose objects are the ordered sets
  $\Delta_n:=\{0<\dots<n\}$ for $n\geq 0$ and whose morphisms are
  non--decreasing functions between them. The category $\psh{\Delta}$ is
  referred to as the category of \emph{simplicial sets}. The canonical inclusion
  $\Delta \hookrightarrow \Cat$ induces the so-called \emph{nerve functor}
  \[
    \begin{aligned}
      N \colon \Cat &\to \psh{\Delta} \\
      \sC &\mapsto \Bigl( \Delta_n \mapsto \Hom_{\Cat}(\Delta_n,\sC) \Bigl).
    \end{aligned}
  \]
  Let us denote by $\cW_{\infty}$ the class of morphisms $u \colon \sA \to \sB$
  of $\Cat$ such that $N(u)$ is a weak homotopy equivalence of simplicial
  sets.\footnote{This means a weak equivalence of the Kan--Quillen model
    structure on simplicial sets.} Recall now the fundamental result of the
  homotopy theory of $\Cat$ \cite[Corollaire 3.3.1]{illusie1972complexe}: the
  nerve functor induces an equivalence at the level of homotopy categories
  \[
    \Ho(\Cat)\simeq\Ho(\psh{\Delta}).
  \]
  In other words, the category $\Cat$ equipped with $\cW_{\infty}$ models
  homotopy types. In the theory of test categories, the point of view is
  reversed and $(\Cat,\cW_{\infty})$ is taken as the fundamental model of
  homotopy types. As it happens, the results of this theory only relies on a few
  formal properties of the class $\cW_{\infty}$, which are shared by other
  classes of weak equivalences (non-necessarily modeling homotopy types),
  referred to as \emph{basic localizers}, and whose definition is recalled
  below.
\end{paragr}
\begin{defn}
  A class $\cW$ of morphisms of $\Cat$ is called a \emph{basic localizer} if it
  satisfies the following properties:
  \begin{rome}
  \item $\cW$ is weakly saturated,
  \item for every small category $\sA$ with a terminal object, the canonical
    morphism to the terminal category
    \[
      \sA \to e
    \]
    is in $\cW$,
  \item for any commutative triangle of $\Cat$,
    \[
      \begin{tikzcd}[column sep=tiny]
        \sA \ar[rr,"u"] \ar[dr] & &\sB \ar[dl]\\
        &\sC,&
      \end{tikzcd}
    \]
    if for every object $c$ of $\sC$, the morphism induced by $u$
    \[
      \sA/c \to \sB/c
    \]
    is in $\cW$, then $u$ is also in $\cW$.
  \end{rome}
\end{defn}
\begin{ex}
  The class $\cW_{\infty}$ is a basic localizer. It is in fact the smallest
  basic localizer \cite{cisinski2004localisateur}. More generally, for any $n
  \geq 0$, let $\cW_n$ be the class of morphisms $u \colon \sA \to \sB$ of
  $\Cat$ such that $N(u)$ induces an equivalences on homotopy groups of
  simplicial sets, up to dimension $n$. Then, $\cW_n$ is a basic localizer
  \cite[Section 9.2]{cisinski2006prefaisceaux} and $(\Cat,\cW_n)$ models
  homotopy $n$\nbd{}types.
\end{ex}
\emph{We now fix once and for all in this section a basic localizer $\cW$ of
  $\Cat$.}
\begin{paragr}
  A small category $\sA$ is called \emph{$\cW$\nbd{}aspherical}, or simply
  \emph{aspherical}, when the canonical morphism to the terminal category $\sA
  \to e$ is in $\cW$.

  More generally, a morphism $u \colon \sA \to \sB$ is
  \emph{$\cW$\nbd{}aspherical}, or simply \emph{aspherical}, if for every object
  $b$ of $\sB$, the category $\sA/b$ is aspherical. Note that it follows from
  the axioms of basic localizers that \emph{every aspherical morphism is in
    $\cW$.} Practically, this is very often how we will prove that a morphism of
  $\Cat$ is a weak equivalence.

  Finally, we say that a small category $\sA$ is \emph{totally
    $\cW$\nbd{}aspherical}, or simply \emph{totally aspherical}, if it is
  aspherical and the diagonal functor $\delta \colon \sA \to \sA \times \sA$ is
  aspherical.
\end{paragr}




For later reference we put here the following lemma. We refer to
\cite[1.1.15]{maltsiniotis2005theorie} for the definition of Grothendieck
fibrations.
\begin{lemma}\label{pullbackasph}
  Let
  \[
    \begin{tikzcd}
      A' \ar[r,"u"] \ar[d,"p'"'] & A \ar[d,"p"] \\
      B' \ar[r,"v"'] & B \ar[from=1-1,to=2-2,"\lrcorner",very near
      start,phantom]
    \end{tikzcd}
  \]
  be a pullback square of $\Cat$. If $p$ is a Grothendieck fibration and $v$ is
  aspherical then $u$ is also aspherical.
\end{lemma}
\begin{proof}
  This is a particular case of the dual of \cite[Théorème
  3.2.15]{maltsiniotis2005theorie}.
\end{proof}

The first step of the theory of test categories is to equip $\psh{\sA}$, for any
small category $\sA$, with a canonical class of weak equivalences. For that, we
use a canonical functor $\psh{\sA} \to \Cat$.
\begin{paragr}
  Let $\sA$ be a small category. We write $i_{\sA} \colon \psh{\sA} \to \Cat$
  for the functor
  \[
    \begin{aligned}
      i_{\sA} \colon \psh{\sA} &\to \Cat \\
      X &\mapsto \sA/X,
    \end{aligned}
  \]
  where $\sA/X$ is the category of elements of $X$. In details, the objects of
  $\sA/X$ are pairs $(a,x)$, where $a$ is an object of $\sA$ and $x \in X(a)$. A
  morphism $(a,x) \to (a',x')$ consists of a morphism $f \colon a \to a'$ in
  $\sA$ such that $x=X(f)(x')$.

  The functor $i_{\sA}$ has a right adjoint given by the following formula
  \[
    \begin{aligned}
      i_{\sA}^* \colon \Cat &\to \psh{\sA}\\
      \sC &\mapsto \Bigl(a \mapsto \Hom_{\Cat}(\sA/a,\sC)\Bigl).
    \end{aligned}
  \]
\end{paragr}
\begin{defn}\label{defn:wepsh}
  Let $\sA$ be a small category. A morphism $\varphi \colon X \to Y$ of
  $\psh{\sA}$ is a \emph{$\cW$\nbd{}equivalence}, or simply a \emph{weak
    equivalence}, if $i_{\sA}(\varphi)$ belongs to $\cW$. We denote by
  $\cW_{\psh{\sA}}$ the class of $\cW$\nbd{}equivalences.

  An object $X$ of $\psh{\sA}$ is \emph{$\cW$\nbd{}aspherical}, or simply
  \emph{aspherical}, if $i_{\sA}(X)$ is a $\cW$\nbd{}aspherical category.

  Finally, an object $X$ of $\psh{\sA}$ is \emph{$\cW$\nbd{}locally aspherical},
  or simply \emph{locally aspherical}, if for every object $a$ of $\sA$, the
  object $X\vert_{\sA/a}$ of $\psh{\sA/a}$, defined as the following
  composition:
  \[
    (\sA/a)^{\op} \to \sA^{\op} \overset{X}{\to} \Set,
  \]
  is aspherical.
\end{defn}
\begin{lemma}
  If $\sA$ is aspherical, then an object $X$ of $\psh{\sA}$ is aspherical if and
  only if the canonical morphism to the terminal presheaf $X \to \ast$ is a weak
  equivalence.
\end{lemma}
\begin{proof}
  See \cite[Section 1.2.6]{maltsiniotis2005theorie}.
\end{proof}
In particular, \emph{when $\sA$ is aspherical}, every locally aspherical
presheaf over $\sA$ is aspherical.
\begin{ex}
  In the case that $\sA=\Delta$, the $\cW_{\infty}$\nbd{}equivalences coincide
  with the usual weak homotopy equivalences, and a simplicial set is
  $\cW_{\infty}$\nbd{}aspherical if and only if it is weakly contractible.
\end{ex}
\begin{defn}[Grothendieck]
  Let $\sA$ be a small category. We say that:
  \begin{enumerate}[label=(\alph*)]
  \item $\sA$ is a \emph{$\cW$\nbd{}pseudo-test category}, or simply a
    \emph{pseudo-test category}, if it satisfies both following conditions:
    \begin{enumerate}[label=(\roman*)]
    \item $\sA$ is aspherical\footnote{In fact, it follows from a result of
        Cisinski \cite[Proposition 4.2.4]{cisinski2006prefaisceaux} and from
        \cite[Proposition 1.3.5]{maltsiniotis2005theorie} that condition (i) is
        implied by condition (ii).},
    \item $i_{\sA} \colon \psh{\sA} \to \Cat$ induces an equivalence at the
      level of homotopy categories
      \[
        \Ho(\psh{\sA}) \simeq \Ho(\Cat),
      \]
    \end{enumerate}
  \item $\sA$ is a \emph{$\cW$\nbd{}weak test category}, or simply a \emph{weak
      test category}, if the adjunction $i_{\sA} \dashv i_{\sA}^*$ is a
    homotopical equivalence (\ref{def:homotopicalequiv}),
  \item $\sA$ is a \emph{$\cW$\nbd{}local test category}, or simply a
    \emph{local test category}, if for every object $a$ of $\sA$, the category
    $\sA/a$ is weak test,
  \item $\sA$ is a \emph{$\cW$\nbd{}test category}, or simply a \emph{test
      category}, if it is both a weak test category and a local test category,
  \item $\sA$ is a \emph{$\cW$\nbd{}strict test category}, or simply a
    \emph{strict test category}, if it is both totally aspherical and a test
    category.
  \end{enumerate}
\end{defn}
\begin{rem}
  We have the following sequence of implications
  \begin{center}
    strict test $\Rightarrow$ test $\Rightarrow$ weak test $\Rightarrow$
    pseudo-test,
  \end{center}
  but it can be shown that the converse of the first two implications do not
  hold. For the converse of the third one, it is still an open question.
\end{rem}
\begin{ex}\label{ex:testcategories}
  The archetypal example of a strict test category is $\Delta$, but it is far
  from being the only one. The class of strict test categories also contains the
  cubical category with connections \cite{maltsiniotis2009categorie}, Joyal's
  $\Theta$ category \cite{cisinski2011categorie}, etc. Examples of test
  categories which are not strict include the cubical category (without
  connections) \cite[Corollaire 8.4.13]{cisinski2006prefaisceaux} and the
  dendroidal category $\Omega$ \cite{ara2019dendroidal}. Examples of weak test
  categories which are not test include the subcategory $\Delta'$ of $\Delta$
  with only monomorphisms as the morphisms \cite[Proposition
  1.7.25]{maltsiniotis2005theorie}. All the examples above are of ``shape-like''
  nature, but it is not always the case. For example, the monoid of
  non-decreasing functions $\bN \to \bN$, seen as a category with only one
  object, is a strict test category \cite[Example
  3.16]{cisinski2011categorie}.
\end{ex}

We now sum up the classical criteria to detect weak test, local test, test and
strict test categories. For details and other characterizations, we refer to the
first chapter of \cite{maltsiniotis2005theorie}.

Recall that we denote by $\Delta_1$ the ordered set $\{0 < 1\}$.
\begin{prop}\label{prop:crittestcat}
  Let $\sA$ be a small category. We have the following characterizations:
  \begin{enumerate}[label={\upshape(\alph*)}]
  \item $\sA$ is a weak test category if and only if for every small category
    $\sC$ with a terminal object, the presheaf $i_{\sA}^*(\sC)$ is aspherical,
  \item $\sA$ is a local test category if and only if the presheaf
    $i_{\sA}^*(\Delta_1)$ is locally aspherical,
  \item $\sA$ is a test category if and only if it is aspherical and
    $i_{\sA}^*(\Delta_1)$ is locally aspherical,
  \item $\sA$ is a strict test category if and only if it is totally aspherical
    and $i_{\sA}^*(\Delta_1)$ is aspherical.
  \end{enumerate}
\end{prop}

Let us end this section with a quick word on aspherical functors and locally
aspherical functors.
\begin{defn}\label{def:asphfun}
  Let $\sA$ be a small category, $i \colon \sA \to \Cat$ a functor and let
  $i^*\colon \Cat \to \psh{\sA}$ be the functor defined as
  \[
    \begin{aligned}
      i^* \colon \Cat &\to \psh{\sA}\\
      \sC &\mapsto \Bigl(a \mapsto \Hom_{\Cat}(i(a),\sC)\Bigr).
    \end{aligned}
  \]
  We say that $i$ is a \emph{$\cW$\nbd{}aspherical functor}, or simply an
  \emph{aspherical functor}, if it satisfies the two following conditions:
  \begin{enumerate}[label={\upshape(\alph*)}]
  \item $i(a)$ has a terminal object for every object $a$ of $\sA$,
  \item if $\sC$ is a small category with a terminal object, $i^*(\sC)$ is an
    aspherical object of $\psh{\sA}$.
  \end{enumerate}

  We say that $i$ is a \emph{$\cW$\nbd{}locally aspherical functor}, or simply a
  \emph{locally aspherical functor}, if it satisfies condition $(a)$ above and
  the following condition instead of $(b)$:
  \begin{itemize}
  \item[$(b')$] $i^*(\Delta_1)$ is a locally aspherical object of $\psh{\sA}$.
  \end{itemize}
\end{defn}
\begin{rem}
  The definition of aspherical functor given above is not the most general
  possible (see \cite[Definition 1.7.1]{maltsiniotis2005theorie}), but it will
  be sufficient for our purpose. (See also \ref{rem:grpdasphfun} below.)
\end{rem}
\begin{rem}\label{rem:localtestfun}
  The first item (resp.\ second item) of \cref{prop:crittestcat} can be
  reformulated as: $\sA$ is a weak test category (resp.\ local test category) if
  and only if $\sA \to \Cat,a \mapsto \sA/a$ is an aspherical functor (resp.\
  locally aspherical functor).
\end{rem}
\begin{prop}
  Let $\sA$ be a small category and $i \colon \sA \to \Cat$ a functor. We have
  the following implications:
  \begin{enumerate}[label={\upshape(\alph*)}]
  \item if $\sA$ is weak test and $i$ is an aspherical functor, then $i^* \colon
    \Cat \to \psh{\sA}$ is a homotopy inverse of $i_{\sA}$,
  \item if $i$ is a locally aspherical functor, then $\sA$ is local test
    category,
  \item if $i$ is a locally aspherical functor and $\sA$ is aspherical, then
    $\sA$ is a test category and $i^*\colon \Cat \to \psh{\sA}$ is a homotopy
    inverse of $i_{\sA}$.
  \end{enumerate}
\end{prop}
\begin{ex}
  The archetypal example of aspherical functor (which is even locally
  aspherical) is the inclusion functor $i \colon \Delta \hookrightarrow \Cat$.
  Then $i^* \colon \Cat \to \psh{\Delta}$ is nothing but the nerve functor.
  Hence, we recover via the above proposition that the nerve induces an
  equivalence of homotopy categories $\Ho(\Cat) \simeq \Ho(\psh{\Delta})$.
\end{ex}
\begin{rem}
  The name \emph{locally aspherical functor} is non standard, but by (b) of the
  previous proposition, it is equivalent to the usual notion of \emph{local test
    functor}. Similarly, an aspherical functor $i \colon \sA \to \Cat$ such that
  $\sA$ is a (weak) test category, is usually called a \emph{(weak) test
    functor}. We will not use this terminology.
\end{rem}

\section{Homotopy theory of presheaves of groupoids}

      \begin{notation}
        We denote by $\Grpd$ the category of (small) groupoids and for a small
        category $\sA$, we denote by $\pgrpd{\sA}$ the category of
        $\Grpd$\nbd{}valued presheaves over $\sA$. That is, $\pgrpd{\sA}$ is the
        category of functors $\sA^{\op} \to \Grpd$ and natural transformations
        between them\footnote{By that, we mean actual \emph{strict} natural
          transformations and not \emph{pseudo} natural transformations.}. We
        use the notation $\ast$ for the terminal object of $\pgrpd{\sA}$.

        The canonical inclusion $\Set \hookrightarrow \Grpd$, which identifies
        sets with discrete groupoids, induces a canonical inclusion $\psh{\sA}
        \hookrightarrow \pgrpd{\sA}$. Hence, every $\Set$\nbd{}valued presheaf
        can be seen as a $\Grpd$\nbd{}valued
        presheaf.
      \end{notation}
      \begin{paragr}\label{paragr:I_A}
        Let $\sA$ be a small category and $X$ an object of $\pgrpd{\sA}$. We
        write $\dbslice{\sA}{X}$ for the \emph{category of elements of $X$},
        which is defined as the following:
        \begin{itemize}[label=-]
        \item an object is a pair $(a,x)$, where $a$ is an object of $\sA$ and
          $x$ is an object of $X(a)$,
        \item a morphism $(a,x) \to (a',x')$ is a pair $(f,k)$, where $f \colon
          a \to a'$ is a morphism of $\sA$, and $k \colon x \overset{\sim}{\to}
          X(f)(x')$ is a morphism of $X(a)$ (which is necessarily an
          isomorphism).
        \end{itemize}
        The identity morphism of $(a,x)$ is given by $(\id_{a},\id_{x})$ and the
        composition of \[(a,x) \overset{(f,k)}{\longrightarrow} (a',x')
          \overset{(f',k')}{\longrightarrow} (a'',x'')\] is given by \[(f'\circ
          f, k'') \colon (a,x) \to (a'',x''),\] where $k''$ is the composite of
        \[
          x \overset{k}{\to} X(f)(x') \overset{X(f)(k')}{\longrightarrow} X(f'
          \circ f)(x'').
        \]
    
        For every object $X$ of $\pgrpd{\sA}$, the category $\dbslice{\sA}{X}$
        comes equipped with a canonical morphism:
        \[
          \begin{aligned}
            \zeta_{X} \colon \dbslice{\sA}{X} &\to \sA\\
            (a,x) &\mapsto a,
          \end{aligned}
        \]
        which is easily checked to be a Grothendieck fibration (see also
        \cref{rem:GrC} below).

        Given a morphism $\alpha \colon X \to X'$ of $\pgrpd{\sA}$ (i.e.\ a
        natural transformation), we define a functor $\dbslice{\sA}{X} \to
        \dbslice{\sA}{X'}$ in the following way:
        \begin{itemize}[label=-]
        \item an object $(a,x)$ of $\dbslice{\sA}{X}$ is sent to the object
          $(a,\alpha_a(x))$ of $X'(a)$,
        \item a morphism $(f,k) \colon (a,x) \to (a',x')$ of $\dbslice{\sA}{X}$
          is sent to the morphism \[(f,\alpha_a(k)) \colon (a,\alpha_a(x)) \to
            (a',\alpha_{a'}(x'))\] of $\dbslice{\sA}{X'}$ (where we used the
          naturality of $\alpha$ for the target of this morphism to be
          compatible).
        \end{itemize}
        This makes the correspondence $X \mapsto \dbslice{\sA}{X}$ functorial in
        $X$, and yields a functor, denoted by $I_{\sA}$:
        \[
          \begin{aligned}
            I_{\sA} \colon \pgrpd{\sA} &\to \Cat \\
            X &\mapsto \dbslice{\sA}{X}.
          \end{aligned}
        \]
        We shall see later that $I_{\sA}$ admits a right adjoint.
      \end{paragr}

      \begin{rem}\label{rem:GrC}
        Via the canonical inclusion $\Grpd \hookrightarrow \Cat$, any
        $\Grpd$\nbd{}valued presheaf $X$ can be seen as a $\Cat$\nbd{}valued
        presheaf. Then, $\dbslice{\sA}{X}$ is nothing but the Grothendieck
        construction of $X$ (see \cref{paragr:GrC} for details) and the
        canonical morphism $\zeta_{X} \colon \dbslice{\sA}{X} \to \sA$ is the
        Grothendieck fibration associated to $X$.
      \end{rem}
      \begin{rem}\label{rem:pshvspgrpd}
        When $X$ is an object of $\psh{\sA}$, which we see as an object of
        $\pgrpd{\sA}$ via the canonical inclusion $\psh{\sA} \hookrightarrow
        \pgrpd{\sA}$, we have
        \[
          i_{\sA}(X)=I_{\sA}(X).
        \]
        In other words, the following triangle is commutative
        \[
          \begin{tikzcd}
            \psh{\sA} \ar[r,hook] \ar[dr,"i_{\sA}"']& \pgrpd{\sA}
            \ar[d,"I_{\sA}"]
            \\
            & \Cat.
          \end{tikzcd}
        \]
      \end{rem}
      \begin{rem}
        By an obvious variation of the Yoneda lemma, the category $I_{\sA}(X)$
        can alternatively be defined as the category whose objects are pairs
        $(a,p \colon a \to X)$ (we identify an object $a$ of $\sA$ with the
        $\Set$\nbd{}valued presheaf represented by $a$) and whose morphisms
        $(a,p) \to (a',p')$ are pairs $(f,\sigma)$, where $f \colon a\to a'$ is
        a morphism of $\sA$ and $\sigma \colon p \overset{\simeq}{\Rightarrow}
        p'\circ f$
        is a natural isomorphism.
      \end{rem}

      \textit{For the rest of this section, we fix once and for all a basic
        localizer $\cW$ of $\Cat$.}
      \begin{defn}
        Let $\sA$ be a small category. A morphism of $\pgrpd{\sA}$
        \[
          \varphi \colon X\to Y
        \]
        is a \emph{$\cW$\nbd{}equivalence}, or simply a \emph{weak equivalence},
        if the induced morphism of $\Cat$
        \[I_{\sA}(\varphi) \colon I_{\sA}(X) \to I_{\sA}(Y)\] is in $\cW$. We
        denote by $\cW_{\pgrpd{\sA}}$ the class of $\cW$\nbd{}equivalences.

        An object $X$ of $\pgrpd{\sA}$ is \emph{$\cW$\nbd{}aspherical}, or
        simply \emph{aspherical}, if the category $I_{\sA}(X)$ is aspherical.
      \end{defn}
      \begin{rem}\label{rem:compatibleasph}
        It follows from \cref{rem:pshvspgrpd} that a $\Set$\nbd{}valued presheaf
        is aspherical in the sense of \cref{defn:wepsh} if only if it is
        aspherical in the sense of the previous definition (using the canonical
        inclusion $\psh{\sA} \hookrightarrow \pgrpd{\sA}$).
      \end{rem}

      In the case that the category $\sA$ is aspherical, there is an equivalent
      characterization of aspherical objects of $\pgrpd{\sA}$ as stated in the
      following lemma.
      \begin{lemma}\label{lemma:asphgprdprsh}
        If $\sA$ is aspherical, then for every object $X$ of $\pgrpd{\sA}$, we
        have the following equivalence
        \begin{center}
          $X \to \termA$ is a weak equivalence $\Leftrightarrow$ $X$ is
          aspherical.
        \end{center}
        Conversely, if this equivalence is true for every $X$ in $\pgrpd{\sA}$,
        then $\sA$ is aspherical.
      \end{lemma}
      \begin{proof}
        Notice $I_{\sA}(\termA)\simeq \sA$, hence $\termA$ is aspherical if and
        only if $\sA$ is aspherical. The equivalence follows then from the
        obvious fact that for a morphism $\varphi \colon X \to Y$ of
        $\pgrpd{\sA}$, if $Y$ is aspherical, then $\varphi$ is a weak
        equivalence if and only if $X$ is aspherical.

        For the second part of the lemma, it suffices to notice that the
        identity morphism $\termA \to \termA$ is always a weak equivalence (as
        all identity morphisms are) and $I_{\sA}(\termA)\simeq \sA$.
      \end{proof}

      For later reference, we put here the following result which gives other
      characterizations of aspherical morphisms of $\Cat$.

\begin{prop}\label{prop:asphmorgrpd}
  Let $u \colon \sA \to \sB$ be a morphism of $\Cat$. The following conditions
  are equivalent:
  \begin{enumerate}[label={\upshape(\alph*)}]
  \item $u$ is aspherical,
  \item the functor $u^* \colon \pgrpd{\sB} \to \pgrpd{\sA}$ preserves and
    reflects aspherical objects, i.e.\ an object $X$ of $\pgrpd{\sB}$ is
    aspherical if and only if $u^*(X)$ is aspherical,
  \item the functor $u^* \colon \pgrpd{\sB} \to \pgrpd{\sA}$ preserves
    aspherical objects, i.e.\ for every aspherical object $X$ of $\pgrpd{\sB}$,
    $u^*(X)$ is aspherical.
  \end{enumerate}
  All these equivalent conditions imply the following condition:
  \begin{enumerate}[label={\upshape(\alph*)},resume]
  \item the functor $u^* \colon \pgrpd{\sB} \to \pgrpd{\sA}$ preserves and
    reflects weak equivalences, i.e.
    \[
      (u^{*})^{-1}(\cW_{\pgrpd{\sA}})=\cW_{\pgrpd{\sB}}.
    \]
  \end{enumerate}
  If moreover $\sA$ and $\sB$ are aspherical, then conditions $(a)$ to $(d)$ are
  all equivalent and equivalent to the following condition:
  \begin{enumerate}[label={\upshape(\alph*)},resume]
  \item the functor $u^*$ preserves weak equivalences, i.e.
    \[
      u^*(\cW_{\pgrpd{\sB}}) \subseteq \cW_{\pgrpd{\sA}}.
    \]
  \end{enumerate}
\end{prop}
\begin{proof}
  Let us start with some preliminaries. It is easily checked that for every
  object $X$ of $\pgrpd{\sB}$, the following square
  \[
    \begin{tikzcd}
      I_{\sA}(u^*(X)) \ar[d,"\zeta_{u^*(X)}"']\ar[r,"\lambda_X"] & I_{\sB}(X)
      \ar[d,"\zeta_{X}"]
      \\
      \sA \ar[r,"u"'] &\sB, \ar[from=1-1,to=2-2,very near start, phantom,
      description,"\lrcorner"]
    \end{tikzcd}
  \]
  where $\lambda_{X}$ is the functor defined on objects as
  \[
    (a,x) \mapsto (u(a),x)
  \]
  and on morphisms as
  \[
    \left( (a,x) \overset{(f,\sigma)}{\longmapsto} (a',x') \right) \mapsto
    \left( (u(a),x) \overset{(u(f),\sigma)}{\longmapsto} (u(a'),x') \right),
  \]
  is a pullback square.

  Now, using the fact the vertical morphisms of the previous pullback square are
  Grothendieck fibrations, it follows from \cref{pullbackasph} that if $u$ is
  aspherical, then $\lambda_{X}$ is aspherical. In particular, with these
  conditions, $X$ is aspherical if and only if $u^*(X)$ is aspherical, which
  proves the implication $(a) \Rightarrow (b)$. The implication $(b) \Rightarrow
  (c)$ is tautological. To prove $(c) \Rightarrow (a)$, consider an object $b$
  of $\sB$, seen as a $\Set$\nbd{}valued representable presheaf (and as an
  object of $\pgrpd{\sA}$ via the inclusion $\psh{\sA} \hookrightarrow
  \pgrpd{\sA}$). It is straightforward to check that $I_{\sA}(u^*(b))\simeq
  \sA/b$. Hence, if condition $(c)$ is satisfied, then $u$ is aspherical.

  \noindent To prove $(c) \Rightarrow (d)$, it suffices to notice that
  $\lambda_X$ is natural in $X$, and thus, for every morphism $f \colon X \to Y$
  of $\pgrpd{\sB}$, we have a commutative square
  \[
    \begin{tikzcd}
      I_{\sA}(u^*(X)) \ar[d,"\lambda_X"'] \ar[r,"I_{\sA}(u^*(f))"] &
      I_{\sA}(u^*(Y)) \ar[d,"\lambda_Y"] \\
      I_{\sB}(X) \ar[r,"I_{\sB}(f)"'] & I_{\sB}(Y).
    \end{tikzcd}
  \]
  We have already seen that if $u$ is aspherical, then the vertical arrows of
  the previous square are weak equivalences. Hence, in this case, $I_{\sB}(f)$
  is a weak equivalence if and only if $I_{\sA}(u^*(f))$ is. By definition of
  weak equivalences of $\Grpd$\nbd{}valued presheaves, this means exactly that
  $u^*$ preserves and reflects weak equivalences.

  The implication $(d) \Rightarrow (e)$ is trivial. Finally, let us prove $(e)
  \Rightarrow (c)$. Let $X$ be an aspherical object of $\pgrpd{\sB}$. Since
  $\sB$ is aspherical, we have already seen (\cref{lemma:asphgprdprsh}) that
  this means exactly that the canonical morphism
  \[
    X \to \ast
  \]
  is a weak equivalence of $\pgrpd{\sB}$. If $u^*$ preserves weak equivalences,
  then
  \[
    u^*(X) \to u^*(\ast)\simeq\ast
  \]
  is a weak equivalence of $\pgrpd{\sA}$. Using that $\sA$ is aspherical, we
  deduce from \cref{lemma:asphgprdprsh} again that $u^*(X)$ is aspherical.
  Hence, $u^*$ preserves aspherical objects.
\end{proof}
\begin{defn}
  A small category $\sA$ is \emph{$\cW$\nbd{}groupoidal pseudo-test}, or simply
  \emph{groupoidal pseudo-test}, if
  \begin{enumerate}[label={\upshape(\alph*)}]
  \item $\sA$ is aspherical,
  \item $I_{\sA}$ induces an equivalence at the level of localized categories:
    \[
      \Ho(\pgrpd{\sA}) \overset{\sim}{\rightarrow} \Ho(\Cat).
    \]
  \end{enumerate}
\end{defn}
\section{Groupoidal weak test categories}
\begin{notation}
  For two (small) categories $\sC$ and $\sD$, we denote by
  $\Homi_{\Cat}(\sC,\sD)$ the groupoid whose objects are functors $\sC \to \sD$
  and whose morphisms are natural \emph{isomorphisms} between those functors.
\end{notation}
\begin{paragr}
  Let $i \colon \sA \to \Cat$ be a functor with $\sA$ a small category. We
  denote by $I^*$ the functor
  \[
    \begin{aligned}
      I^* \colon \Cat &\to \pgrpd{\sA}\\
      \sC &\mapsto \Homi_{\Cat}(i(-),\sC).
    \end{aligned}
  \]
  When $i$ is the functor $\sA \to \Cat, a \mapsto \sA/a$, we use the special
  notation $I_{\sA}^*$ for the functor $I^*$. In other words, for a small
  category $\sC$, $I_{\sA}^*$ is the functor
  \[
    \begin{aligned}
      I_{\sA}^* \colon \Cat &\to \pgrpd{\sA}\\
      \sC &\mapsto \Homi_{\Cat}(\sA/(-),\sC).
    \end{aligned}
  \]
\end{paragr}
The following lemma is without a doubt a folkloric result and I claim no
originality for it.
\begin{lemma}\label{lemma:I*adjoint}
  Let $i \colon \sA \to \Cat$ be a functor, with $\sA$ a small category. The
  functor $I^* \colon \Cat \to \pgrpd{\sA}$ has a left adjoint. Moreover, in the
  case that $i$ is the functor $a \mapsto \sA/a$ (and thus $I^*=I_{\sA}^*$),
  this left adjoint is $I_{\sA} \colon \pgrpd{\sA} \to \Cat$
  (\cref{paragr:I_A}).
\end{lemma}
\begin{proof}
  Let us denote by $\int^{a \in \sA}F(a,a)$ (resp.\ $\int_{a \in \sA}F(a,a)$)
  the co-end (resp.\ the end) of a functor $F \colon \sA^{\op}\times \sA \to
  \Cat$. We define a functor $I_! \colon \pgrpd{\sA} \to \Cat$ as
  \[
    X \mapsto \int^{a \in \sA}X(a)\times i(a).
  \]

  For every small category $\sC$ and every $\Grpd$\nbd{}valued presheaf $X$, we
  then have the following sequence of natural isomorphisms:
  \[
    \begin{aligned}
      \Hom_{\Cat}(\int^{a \in \sA}X(a)\times i(a),\sC) &\simeq \int_{a \in
        \sA}\Hom_{\Cat}(X(a)\times i(a),\sC)\\
      &\simeq \int_{a \in \sA}\Hom_{\Grpd}(X(a),\Homi_{\Cat}(i(a),\sC)) \\
      &\simeq \Hom_{\pgrpd{\sA}}(X,I^*(\sC)).
    \end{aligned}
  \]
  Hence, $I_!$ is left adjoint of $I^*$. In the case that $i \colon \sA \to
  \Cat$ is the functor $a \mapsto \sA/a$, we have
  \[
    I_!(X)=\int^{a \in \sA}X(a)\times \sA/a,
  \]
  which is nothing but the Grothendieck construction of $X$. To see this, recall
  that the Grothendieck construction of a $\Cat$\nbd{}valued functor (and,
  \textit{a fortiori}, for a $\Grpd$\nbd{}valued functor) is its oplax colimit
  \cite{gray1969categorical}, and that the oplax colimit of a contravariant
  functor is computed as the colimit weighted by the slices of the source
  \cite{street1976limits}. The conclusion follows then from \cref{rem:GrC}.
\end{proof}
\textit{We now fix, once and for all in this section, a basic localizer $\cW$ of
  $\Cat$.}

\begin{defn}\label{def:grpdweaktest}
  A small category $\sA$ is \emph{$\cW$\nbd{}groupoidal weak test}, or simply
  \emph{groupoidal weak test}, if the adjunction $I_{\sA} \dashv I_{\sA}^*$ is a
  homotopical equivalence (\cref{def:homotopicalequiv}).
\end{defn}
\begin{rem}
  An immediate computation shows that $I_{\sA}I_{\sA}^*(e) \simeq \sA$. Thus, if
  $\sA$ is groupoidal weak test, the co-unit morphism $\sA \to e$ is a weak
  equivalence and so $\sA$ is aspherical. This shows that every groupoidal weak
  test category is a groupoidal pseudo-test category.
\end{rem}

We would like now to find characterizations of groupoidal weak test categories.
For that, we begin by studying a class of homotopically well-behaved functors
$\sA \to \Cat$.
\begin{defn}\label{def:grpdasphfun}
  Let $\sA$ be a small category. A functor $i \colon \sA \to \Cat$ is
  \emph{$\cW$\nbd{}groupoidal aspherical}, or simply \emph{groupoidal
    aspherical}, if:
  \begin{enumerate}[label={\upshape(\alph*)}]
  \item for every object $a$ of $\sA$, the category $i(a)$ has a terminal
    object,
  \item for every small category $\sC$ with a terminal object, the
    $\Grpd$\nbd{}valued presheaf $I^*(\sC)$ is aspherical.
  \end{enumerate}
\end{defn}
\begin{rem}\label{rem:grpdasphfun}
  A more general notion of groupoidal aspherical functor is obtained by
  replacing the conditions of the previous definitions by:
  \begin{enumerate}[label=(\alph*')]
  \item $i(a)$ is aspherical for every object $a$ of $\sA$,
  \item for every small aspherical category $\sC$, $I^*(\sC)$ is aspherical,
  \end{enumerate}
  which is the straightforward generalization of \cite[Definition
  1.7.1]{maltsiniotis2005theorie}. (We will see in \cref{prop:critgrpdasphfun}
  that when (a) is satisfied, condition (b) and (b') are equivalent). However,
  the author of the notes does not know if the theory fully works for this more
  general notion of groupoidal aspherical functor, and the restricted version we
  chose is sufficient for our purpose. For the interested reader, the difficulty
  is that it seems that \cite[Lemme 1.7.4 (b)]{maltsiniotis2005theorie} cannot
  by generalized to $\Grpd$\nbd{}valued presheaves. (Nevertheless, see
  \cref{lemma174grpd}(b) below for a \emph{partial} generalization of this
  lemma.)
\end{rem}

\begin{paragr}\label{yogatermobj}
  Let $i \colon \sA \to \Cat$ be a functor, with $\sA$ a small category, and
  suppose that for every object $a$ of $\sA$, the category $i(a)$ has a terminal
  object $e_{a}$. We are going to define a canonical natural transformation
  \[
    \alpha \colon I_{\sA}I^* \Rightarrow \id_{\Cat}.
  \]
  Let $\sC$ be a small category. Spelling out the definitions, we see that the
  category $I_{\sA}I^*(\sC)$ has for object pairs $(a,p \colon i(a) \to \sC)$
  where $a$ is an object of $\sA$ and $p$ a morphism of $\Cat$. And a morphism
  $(a,p) \to (a',p')$ in $I_{\sA}I^*(\sC)$ consists of a pair $(f,\sigma)$,
  where $f \colon a \to a'$ is a morphism of $\sA$ and
  \[
    \sigma \colon p \overset{\simeq}{\Rightarrow} p' \circ i(f)
  \]
  is a natural isomorphism.

  At the level of objects, we define the morphism $\alpha_{\sC} \colon
  I_{\sA}I^*(\sC) \to \sC$ with the formula
  \[
    \alpha_{\sC}(a,p):=p(e_a).
  \]
  At the level of morphisms, the image of a morphism $(f,\sigma) \colon (a,p)
  \to (a',p')$ is defined as the composite
  \[
    \alpha_{\sC}(f,\sigma):=p(e_a) \overset{\sigma_{e_a}}{\longrightarrow}
    p'(i(f)(e_a)) \to p'(e_{a'}),
  \]
  where the morphism on the right is induced by the canonical morphism
  $i(f)(e_a) \to e_{a'}$.

  We leave it to the reader to check that $\alpha_{\sC} \colon I_{\sA}I^*(\sC)
  \to \sC$ does indeed define a functor and that it is natural in $\sC$.
\end{paragr}
\begin{rem}\label{alphaiscounit}
  When $i \colon \sA \to \Cat$ is the functor $a \mapsto A/a$, then $\alpha
  \colon I_{\sA}I^*_{\sA} \Rightarrow \id_{\Cat}$ is nothing but the co--unit of
  the adjunction $I_{\sA} \dashv I_{\sA}^*$.
\end{rem}
\begin{lemma}\label{lemma:slicecommut}
  Let $i \colon A \to \Cat$ be a functor with $\sA$ a small category and suppose
  that for every object $a$ of $\sA$, the category $i(a)$ has a terminal object.
  For every small category $\sC$ and $c$ an object of $\sC$, we have a canonical
  isomorphism
  \[
    I_{\sA}I^*(\sC)/c\simeq I_{\sA}I^*(\sC/c).
  \]
\end{lemma}
\begin{proof}
  Let us denote by $e_a$ the terminal object of $i(a)$. The category
  $I_{\sA}I^*(\sC)/c$ is described as follows:
  \begin{itemize}[label=-]
  \item an object is a triple $(a,p \colon i(a) \to \sC, g \colon p(e_a) \to
    c)$, where $a$ is an object of $\sA$, $p$ is a morphism of $\Cat$, and $g$
    is a morphism of $\sC$,
  \item a morphism $(a,p,g) \to (a',p',g')$ is a couple $(f, \sigma)$, where $f
    \colon a \to a'$ is a morphism of $\sA$ and
    \[
      \begin{tikzcd}
        i(a) \ar[r,"i(f)"] \ar[dr,"p"',""{name=A,right}] & i(a') \ar[d,"p'"] \\
        & \sC \ar[from=A,to=1-2,"\sigma","\simeq"',Rightarrow]
      \end{tikzcd}
    \]
    is a natural isomorphism, such that the triangle
    \[
      \begin{tikzcd}
        p(e_a) \ar[r,"{\alpha_{\sC}(f,\sigma)}"] \ar[dr,"g"'] & p'(e_a')
        \ar[d,"g'"]\\
        &c
      \end{tikzcd}
    \]
    is commutative.
  \end{itemize}
  The category $I_{\sA}I^*(\sC/c)$ is described as:
  \begin{itemize}[label=-]
  \item an object is a couple $(a,q \colon i(a) \to \sC/c)$, where $a$ is an
    object of $\sA$ and $q$ is a morphism of $\Cat$,
  \item a morphism $(a,q) \to (a',q')$ is a pair $(f,\sigma)$, where $f \colon a
    \to a'$ is a morphism of $\sA$ and
    \[
      \begin{tikzcd}
        i(a) \ar[r,"i(f)"] \ar[dr,"q"',""{name=A,right}] & i(a') \ar[d,"q'"] \\
        & \sC/c \ar[from=A,to=1-2,"\sigma","\simeq"',Rightarrow]
      \end{tikzcd}
    \]
    a natural isomorphism.
  \end{itemize}
  Let us write $\pi_c \colon \sC/c \to \sC$ for the canonical projection
  functor. We define a morphism of $\Cat$ as,
  \[
    \begin{aligned}
      \theta \colon I_{\sA}I^*(\sC/c) &\to I_{\sA}I^*(\sC)/c\\
      (a,q \colon i(a) \to \sC/c) &\mapsto (a,\pi_c\circ q \colon i(a) \to \sC,
      q(e_a) \colon \pi_c(q(e_a)) \to c),
    \end{aligned}
  \]
  the definition on morphisms being the obvious one.
  
  We now leave to the reader to verify that this morphism is indeed an
  isomorphism. This mainly amounts to showing the following general fact: let
  $\sA$ and $\sB$ be two categories and suppose that $\sA$ has a terminal object
  $t_{\sA}$. For any object $b$ of $\sB$, a functor $f \colon \sA \to \sB/b$ is
  entirely determined by the post-composition $\sA \overset{f}{\to} \sB/b \to
  \sB$ and by $f(t_{\sA})$, seen as a morphism of $\sB$ whose target is $b$.
\end{proof}
\begin{prop}\label{prop:critgrpdasphfun}
  Let $i \colon \sA \to \Cat$ be a functor such that for every $a$ in $\sA$ the
  category $i(a)$ has a terminal object. The following conditions are
  equivalent:
  \begin{enumerate}[label={\upshape(\alph*)}]
  \item $i$ is groupoidal aspherical,
  \item $I^*$ preserves aspherical objects, i.e.\ for every small aspherical
    category $\sC$, the $\Grpd$\nbd{}valued presheaf $I^*(\sC)$ is aspherical,
  \item $I^*$ preserves and reflects aspherical objects, i.e.\ a small category
    $\sC$ is aspherical if and only if the $\Grpd$\nbd{}valued presheaf
    $I^*(\sC)$ is aspherical,
  \item for every small category $\sC$, the canonical morphism
    \[
      \alpha_{\sC} \colon I_{\sA}I^*(\sC) \to \sC
    \]
    is a weak equivalence,
  \item $\sA$ is aspherical and $I^*$ preserves and reflects weak equivalences,
    i.e.\
    \[
      (I^{*})^{-1}(\cW_{\pgrpd{\sA}})=\cW,
    \]
  \item $\sA$ is aspherical and $I^*$ preserves weak equivalences, i.e.\
    \[
      I^*(\cW)\subseteq \cW_{\pgrpd{\sA}}.
    \]
  \end{enumerate}
\end{prop}
\begin{proof}
  The implications $(d) \Rightarrow (c) \Rightarrow (b) \Rightarrow (a)$ are
  trivial (the last one comes from that every category with a terminal object is
  aspherical). For the implication $(a) \Rightarrow (d)$, it follows from
  \cref{lemma:slicecommut} (and the fact the slices $\sC/c$ have a terminal
  object) that $\alpha_{\sC}$ is aspherical hence a weak equivalence. This
  proves the equivalence of the first four conditions. For the implication $(d)
  \Rightarrow (e)$, notice first that if $\sC=e$ is the terminal category, we
  have $I_{\sA}I^*(e)\simeq \sA$, and hence $(c)$ implies that $\sA$ is
  aspherical. The fact that $(d)$ implies that $I^*$ preserves and reflects weak
  equivalences follows from the naturality of $\alpha$, 2-out-of-3 and the fact
  that $\cW_{\pgrpd{\sA}}=I_{\sA}^{-1}(\cW)$ by definition. The implication $(e)
  \Rightarrow (f)$ is tautological. Finally, for the implication $(f)
  \Rightarrow (b)$, let $\sC$ be a small aspherical category and consider the
  canonical morphism $\sC \to e$, which is, by definition, a weak equivalence.
  Since $I^*$ preserves weak equivalences by hypothesis, the induced morphism
  \[
    I^*(\sC) \to I^*(e)
  \]
  is a weak equivalence of $\pgrpd{\sA}$. But $I_{\sA}I^*(e)\simeq \sA$, which
  is by hypothesis aspherical. Hence, $I^*(e)$ is aspherical and then so is
  $I^*(\sC)$.
\end{proof}
We then obtain the following characterization of groupoidal weak test
categories.
\begin{prop}\label{prop:critgrpdweaktest}
  Let $\sA$ be a small category. The following conditions are equivalent:
  \begin{enumerate}[label={\upshape(\alph*)}]
  \item $\sA$ is groupoidal weak test,
  \item for every small category $\sC$ with a terminal object, the
    $\Grpd$\nbd{}valued presheaf $I_{\sA}^*(\sC)$ is aspherical.
  \end{enumerate}
\end{prop}
\begin{proof}
  The second condition means exactly that the functor $\sA \to \Cat, a \mapsto
  \sA/a$ is groupoidal aspherical (since the slice categories $\sA/a$ have a
  terminal object). The equivalence follows then from condition (d) of
  \cref{prop:critgrpdasphfun} combined with \cref{alphaiscounit} and
  \cref{lemma:crithomotopicalequiv}.
\end{proof}
Interestingly, we also obtain the following result.
\begin{cor}\label{cor:grpdasphfunhmtpyinverse}
  Let $\sA$ be a groupoidal weak test category and $i \colon \sA \to \Cat$ a
  groupoidal aspherical functor. Then, $I^* \colon \Cat \to \pgrpd{\sA}$ is a
  homotopy inverse of $I_{\sA} \colon \pgrpd{\sA} \to \Cat$.
\end{cor}
\begin{proof}
  We already know that if $i$ is a groupoidal aspherical functor, we have a
  natural transformation $\alpha \colon I_{\sA}I^* \Rightarrow \id_{\Cat}$ which
  is a weak equivalence argument by argument. Consider now the following zigzag
  of natural transformations
  \[
    I^*I_{\sA} \Rightarrow I_{\sA}^*I_{\sA}I^*I_{\sA} \Rightarrow
    I_{\sA}^*I_{\sA} \Leftarrow \id,
  \]
  where the arrow on the left is induced by the unit of the adjunction $I_{\sA}
  \dashv I_{\sA}^*$, the middle arrow is obtained by post-composing $I_{\sA}^*$
  to $\alpha$ and pre-composing with $I_{\sA}$ and the arrow on the right is
  again the unit of the adjunction $I_{\sA} \dashv I_{\sA}^*$. Since $\sA$ is
  groupoidal weak test, the unit of this adjunction is a weak equivalence
  argument by argument and $I_{\sA}^*$ preserves weak equivalences. This proves
  that the three natural transformations of the previous zigzag are weak
  equivalences argument by argument.
\end{proof}
\begin{rem}
  Following the terminology from the usual test category theory, a groupoidal
  aspherical functor whose source is a weak test category ought to be called a
  \emph{groupoidal weak test functor}.
\end{rem}
\begin{rem}
  Even if we remove the hypothesis that $\sA$ is groupoidal weak test from the
  previous corollary, it follows trivially from \cref{prop:critgrpdasphfun}(d)
  that $I^*$ is a homotopical inverse ``on one side'' of $I_{\sA}$, but it does
  not seem to be a homotopical inverse ``on both sides'' in general. No
  counter-example, however, is known by the author of this paper.
\end{rem}
\section{Groupoidal test categories}
\emph{We fix once and for all in this section a basic localizer $\cW$ of
  $\Cat$.}
\begin{defn}
  A small category $\sA$ is \emph{$\cW$-groupoidal local test}, or simply
  \emph{groupoidal local test}, if for every object $a$ of $\sA$, the category
  $\sA/a$ is groupoidal weak test. We say that $\sA$ is \emph{$\cW$-groupoidal
    test}, or simply \emph{groupoidal test}, if it is both a groupoidal weak
  test category and a groupoidal local test category.
\end{defn}
We are now going to look for characterizations of groupoidal local test
categories and groupoidal test categories. For that, it is useful to introduce
first a variation of the notion of aspherical object of $\pgrpd{\sA}$.
\begin{notation}
  Let $\sA$ be a small category and $a$ an object of $\sA$. Given a
  $\Grpd$\nbd{}valued presheaf $X$ over $A$, we denote by $X\vert_{\sA/a}$ the
  $\Grpd$\nbd{}valued presheaf over $\sA/a$ defined as the composition
  \[
    (\sA/a)^{\op} \to \sA^{\op} \overset{X}{\to} \Grpd,
  \]
  where $\sA/a \to \sA$ is the canonical projection.
\end{notation}

\begin{defn}\label{def:grpdlocallyasph}
  Let $\sA$ be a small category. An object $X$ of $\pgrpd{\sA}$ is
  \emph{$\cW$\nbd{}locally aspherical}, or simply \emph{locally aspherical}, if
  for every object $a$ of $\sA$, the object $X\vert_{\sA/a}$ of
  $\pgrpd{(\sA/a)}$ is aspherical.
\end{defn}

We now have the following reformulation.
\begin{prop}\label{prop:basiccritgrpdlocaltest}
  Let $\sA$ be a small category. The following conditions are equivalent:
  \begin{enumerate}[label={\upshape(\alph*)}]
  \item $\sA$ is groupoidal local test,
  \item for every small category $\sC$ with a terminal object, $I_{\sA}^*(\sC)$
    is locally aspherical.
  \end{enumerate}
\end{prop}
\begin{proof}
  Follows immediately from \cref{prop:critgrpdweaktest} and the fact that for
  every small category $\sC$ and every object $a$ of $\sA$, we have
  \[
    I_{\sA/a}^*(\sC) \simeq I_{\sA}^*(\sC)\vert_{\sA/a}.\qedhere
  \]
\end{proof}
More generally, we can consider the following variation of groupoidal aspherical
functor.
\begin{defn}
  Let $\sA$ be a small category. A functor $i \colon \sA \to \Cat$ is
  \emph{\cW-groupoidal locally aspherical}, or simply \emph{groupoidal locally
    aspherical}, if the following conditions are satisfied:
  \begin{enumerate}[label={\upshape(\alph*)}]
  \item for every object $a$ of $\sA$, $i(a)$ has a terminal object,
  \item for every small category $\sC$ with a terminal object, $I^*(\sC)$ is
    locally aspherical.
  \end{enumerate}
\end{defn}

\begin{rem}\label{rem:grpdlocaltestfun}In other words, \cref{prop:basiccritgrpdlocaltest} says that $\sA$ is
  groupoidal local test if and only if $\sA \to \Cat, a \mapsto \sA/a$ is a
  groupoidal locally aspherical functor.
\end{rem}

We now turn to a couple of key technical results on locally aspherical
$\Grpd$\nbd{}valued presheaves.
   
\begin{lemma}\label{lemma:isolocasph}
  Let $\sA$ be a small category, $a$ an object of $\sA$ and $X$ an object of
  $\pgrpd{\sA}$. We have a sequence of isomorphisms natural in $X$
  \[
    I_{\sA}(a\times X) \simeq I_{\sA/a}(X\vert_{\sA/a}) \simeq I_{\sA}(X)/a,
  \]
  where:
  \begin{itemize}[label=-]
  \item on the left hand side, we abusively wrote $a$ for the $\Set$-valued
    presheaf represented by $a$,
  \item on the right hand side, we implicitly used the canonical morphism
    $\zeta_X \colon I_{\sA}(X) \to \sA$ to make sense of the slice category.
  \end{itemize}
\end{lemma}
\begin{proof}
  These three categories have the same description (up to canonical
  isomorphism):
  \begin{itemize}[label=-]
  \item an object is a triple $(a',p \colon a' \to a, x)$, where $a'$ is an
    object of $\sA$, $p$ is a morphism of $\sA$, and $x$ is an object of
    $X(a')$,
  \item a morphism $(a',p,x) \to (a'',p',x')$ is a pair $(f,k)$ where $f \colon
    a' \to a''$ is a morphism of $\sA$ such that $p'\circ f= p$ and $k \colon x
    \to X(f)(x')$ is a morphism of $X(a')$.\qedhere
  \end{itemize}
\end{proof}
\begin{lemma}\label{lemma:characgrpdlocasph}
  Let $\sA$ be a small category and $X$ an object of $\pgrpd{\sA}$. The
  following conditions are equivalent:
  \begin{enumerate}[label={\upshape(\alph*)}]
  \item $X$ is locally aspherical,
  \item the image by $I_{\sA}$ of the morphism to the terminal object $X \to
    \termA$,
    \[
      I_{\sA}(X)\to \sA,
    \]
    is an aspherical morphism of $\Cat$,
  \item for every object $a$ of $\sA$, seen as a representable presheaf (and as
    an object of $\pgrpd{\sA}$, via the canonical inclusion $\psh{\sA}
    \hookrightarrow \pgrpd{\sA}$), the product in $\pgrpd{\sA}$
    \[
      a \times X
    \]
    is aspherical,
  \item the morphism $X \to \ast$ is a universal weak equivalence, i.e.\ for
    every object $Y$ of $\pgrpd{\sA}$, the canonical projection
    \[
      Y\times X \to Y
    \]
    is a weak equivalence.
  \end{enumerate}
\end{lemma}
\begin{proof}
  The equivalence of $(a)$, $(b)$ and $(c)$ follows immediately from
  \cref{lemma:isolocasph}. For the implication $(b) \Rightarrow (d)$, consider
  the following cartesian square
  \[
    \begin{tikzcd}
      Y\times X \ar[r] \ar[d] & X \ar[d] \\
      Y \ar[r]& \ast. \ar[from=1-1,to=2-2,phantom,"\lrcorner",very near start]
    \end{tikzcd}
  \]
  We leave it to the reader to check that the functor $I_{\sA}$ preserves
  cartesian squares. Hence, we obtain a cartesian square
  \[
    \begin{tikzcd}
      I_{\sA}(Y\times X) \ar[r] \ar[d] & I_{\sA}(X) \ar[d,"\zeta_X"] \\
      I_{\sA}(Y) \ar[r,"\zeta_{Y}"']& \sA.
      \ar[from=1-1,to=2-2,phantom,"\lrcorner", very near start]
    \end{tikzcd}
  \]
  Because $I_{\sA}(Y) \to \sA$ is a Grothendieck fibration, we deduce from
  \cref{pullbackasph} that the left vertical morphism of the previous square is
  an aspherical morphism of $\Cat$, which implies that $Y\times X \to X$ is a
  weak equivalence. Finally, for the implication $(d) \Rightarrow (c)$, let $a$
  be an object of $\sA$, and consider the canonical projection
  \[
    a\times X \to a
  \]
  (where once again we wrote $a$ for the $\Set$-valued presheaf represented by
  $a$). Since the presheaf $a$ is aspherical (because the category
  $I_{\sA}(a)=\sA/a$ has a terminal object), it follows that $a\times X$ is
  aspherical.
\end{proof}
\begin{rem}
  \cref{lemma:characgrpdlocasph} hides a subtlety between the classical theory
  for $\Set$\nbd{}valued presheaves and the theory for $\Grpd$\nbd{}valued
  presheaves. Indeed, let us call \emph{local weak equivalence} a morphism $f
  \colon X \to Y$ of $\pgrpd{\sA}$ such that for every object $a$ of $\sA$, the
  morphism of $\pgrpd{(\sA/a)}$
  \[
    f\vert_{\sA/a} \colon X\vert_{\sA/a} \to Y\vert_{\sA/a}
  \]
  is a weak equivalence. Then, an object $X$ of $\pgrpd{\sA}$ is locally
  aspherical if and only if $X \to \ast$ is a local weak equivalence. Now,
  \cref{lemma:characgrpdlocasph}(d) tells us that $X$ is locally aspherical if
  and only if $X \to \ast$ is a universal weak equivalence, and we might think
  that this characterization is true for all local weak equivalences (as the
  analogue result is true for $\Set$\nbd{}valued presheaves \cite[Proposition
  1.2.5]{maltsiniotis2005theorie}). This does not work though. The correct
  generalization, which goes beyond the scope of this paper, involves what ought
  to be called ``comma-universal weak equivalence'', whose definition is the
  same as universal equivalence only pullback squares are replaced by a comma
  squares.
\end{rem}

     \begin{lemma}\label{cor:locasphisasph}
       Let $\sA$ be a small category and $X$ an object of $\pgrpd{\sA}$. If
       $\sA$ is aspherical, then we have the following implication
       \begin{center}
         $X$ locally aspherical $\Rightarrow$ $X$ aspherical.
       \end{center}
     \end{lemma}
      
      \begin{proof}
        Thanks to \cref{lemma:characgrpdlocasph} (for example item $(b)$), we
        know that if an object $X$ of $\pgrpd{\sA}$ is locally aspherical, then
        $X \to \termA$ is a weak equivalence. The conclusion follows from
        \cref{lemma:asphgprdprsh}.
      \end{proof}
      From this last lemma, we immediately deduce the following two
      propositions.
      \begin{prop}
        Let $\sA$ be a small category and $i \colon \sA \to \Cat$ a groupoidal
        locally aspherical functor. Then $i$ is a groupoidal aspherical functor
        if and only if $\sA$ is aspherical.
      \end{prop}
      \begin{prop}\label{prop:localtestandasph}
        Let $\sA$ be a groupoidal local test category. Then $\sA$ is groupoidal
        test if and only if it is aspherical.
      \end{prop}

      This last result means that in order to characterize groupoidal test
      categories, it suffices to characterize groupoidal local test categories.
   
      \begin{paragr}
        Let $\sM$ be a category with finite products (this includes a terminal
        object, which we denote by $e_{\sM}$). An \emph{interval} in $\sM$ is a
        triple $(\bI,i_0,i_1)$, where $\bI$ is an object of $\sM$, and $i_0$ and
        $i_1$ are morphisms of $\sM$ from $e_{\sM}$ to $\bI$
        \[
          i_0,i_1 \colon e_{\sM} \to \bI.
        \]
      
        A \emph{morphism of intervals} $(\bI,i_0,i_1) \rightarrow
        (\bI',i_0',i_1')$ consists of a morphism $\varphi \colon \bI \to \bI'$
        of $\sM$ such that $i_{\varepsilon}'=\varphi\circ i_{\varepsilon}$ for
        $\varepsilon=0,1$.

        Now, let $f,g \colon X \to Y$ be two parallel morphisms of $\sM$ and
        $(\bI,i_0,i_1)$ an interval in $\sM$. A \emph{$\bI$\nbd{}homotopy from
          $f$ to $g$} is a morphism $h \colon \bI\times X \to Y$ of $\sM$ such
        that the following diagram is commutative
        \[
          \begin{tikzcd}
            X \ar[dr,"f"] \ar[d,"i_0\times X"'] & \\
            \bI\times X \ar[r,"h"]& Y.\\
            X \ar[ru,"g"'] \ar[u,"i_1\times X"] &
          \end{tikzcd}
        \]
        We consider the smallest equivalence relation on the set
        $\Hom_{\sM}(X,Y)$ such that $f$ is equivalent to $g$ if there exists a
        $\bI$\nbd{}homotopy from $f$ to $g$. If two morphisms $X \to Y$ are in
        the same equivalence class for this relation, we say that they are
        \emph{$\bI$\nbd{}homotopic}.

        We say that an object $X$ of $\sM$ is \emph{$\bI$\nbd{}contractible} if
        $\id_X \colon X \to X$ is $\bI$\nbd{}homotopic to a constant morphism
        (i.e.\ a morphism which factorizes through the terminal object
        $e_{\sM}$).

        Finally, notice that if $F \colon \sM \to \sM'$ is a functor preserving
        finite products, then for any interval $(\bI,i_0,i_1)$ of $\sM$,
        $(F(\bI),F(i_0),F(i_1))$ is a interval of $\sM'$, and $F$ sends
        $\bI$\nbd{}homotopic morphisms to $F(\bI)$\nbd{}homotopic morphisms. In
        particular, $F$ sends $\bI$\nbd{}contractible objects to
        $F(\bI)$\nbd{}contractible objects.
      \end{paragr}
      \begin{ex}
        Let $\Delta_1$ be the poset $\{0 < 1\}$ seen as an object of $\Cat$, and
        denote by $e_0,e_1 \colon e \to \Delta_1$ the canonical inclusions of
        $0$ and $1$ respectively. Then $(\Delta_1,e_0,e_1)$ is an interval of
        $\Cat$. A $\Delta_1$\nbd{}homotopy from a morphism $u \colon A \to B$ to
        a morphism $v \colon A \to B$ is nothing but a natural transformation $u
        \Rightarrow v$. Notice that a small category with either a terminal
        object or an initial object is $\Delta_1$\nbd{}contractible.
      \end{ex}
      \begin{ex}
        Let $i \colon \sA \to \Cat$ be a functor where $\sA$ is a small
        category. Since $I^*$ preserves limits,
        $(I^*(\Delta_1),I^*(e_0),I^*(e_1))$ is an interval of $\pgrpd{\sA}$.
      \end{ex}
      The following lemma relates the notion of $\bI$\nbd{}homotopy with the
      homotopy theory induced by a class of weak equivalences in the ambient
      category.
      \begin{lemma}\label{lemma:homotopy}
        Let $\sM$ be a category with finite products, $\sW$ a weakly saturated
        class of morphisms of $\sM$ and $\bI$ an interval of $\sM$ such that the
        canonical morphism to the terminal object $\bI \to e_{\sM}$ is
        universally in $\sW$. Then, for every $\bI$\nbd{}contractible object $X$
        of $\sM$, the canonical morphism $X \to e_{\sM}$ is universally in
        $\sW$.
      \end{lemma}
      \begin{proof}
        This is a reformulation of \cite[1.4.6]{maltsiniotis2005theorie}.
      \end{proof}
      We can now apply this to groupoidal locally aspherical functors.
      \begin{prop}\label{prop:criteriongrpdasphfun}
        Let $i \colon \sA \to \Cat$ be a functor, with $\sA$ a small category,
        such that for every object $a$ of $\sA$, the category $i(a)$ has a
        terminal object. The following conditions are equivalent:
        \begin{enumerate}[label={\upshape(\alph*)}]
        \item $i$ is groupoidal locally aspherical,
        \item $I^*(\Delta_1)$ is locally aspherical.
        \end{enumerate}
      \end{prop}
      \begin{proof}
        The implication $(a) \Rightarrow (b)$ is trivial because $\Delta_1$ has
        terminal object. For the converse, we need to show that for every small
        category $\sC$ with a terminal object, $I^*(\sC)$ is locally aspherical.
        By \cref{lemma:characgrpdlocasph}(d), this is equivalent to showing that
        $I^*(\sC)\to \ast$ is a universal weak equivalence. Notice that $I^*$
        preserves limits and so it sends $\Delta_1$\nbd{}contractible objects of
        $\Cat$ to $I^*(\Delta_1)$\nbd{}contractible objects of $\pgrpd{\sA}$.
        Since every small category with a terminal object is
        $\Delta_1$\nbd{}contractible, the result follows from
        \cref{lemma:homotopy}.
      \end{proof}

      We could now apply the previous proposition to the functor $\sA \to \Cat,
      a \mapsto \sA/a$ and obtain a characterization of groupoidal local test
      categories. As it happens, we will soon obtain an even finer
      characterization, but we first need some more results on intervals.

      \begin{defn}
        Let $\sM$ be a category with finite products, and whose terminal object
        is denoted by $e_{\sM}$. A \emph{multiplicative interval} in $\sM$ is a
        interval $(\bL,\lambda_0,\lambda_1)$ together with a binary operation
        \[
          \Lambda \colon \bL\times \bL \to \bL,
        \]
        such that $\lambda_0$ is a unit on the left and $\lambda_1$ is absorbing
        on the left. In other words, the following two diagrams are commutative:
        \[
          \begin{tikzcd}[column sep=huge]
            e_{\sM}\times \bL \ar[r,"\lambda_0\times \id_\bL"]
            \ar[dr,"\id_\bL"'] & \bL\times
            \bL \ar[d,"\Lambda"] \\
            & \bL,
          \end{tikzcd}
          \qquad
          \begin{tikzcd}[column sep=huge]
            e_{\sM} \times \bL \ar[r,"\lambda_1\times \id_\bL"] \ar[d] &
            \bL\times \bL
            \ar[d,"\Lambda"] \\
            e_{\sM} \ar[r,"\lambda_1"'] & \bL.
          \end{tikzcd}
        \]
      \end{defn}
      \begin{ex}\label{ex:multinterval}
        The interval ($\Delta_1,e_0,e_1$) of $\Cat$ is multiplicative when
        equipped with the binary operation
        \[
          \begin{aligned}
            \Delta_1\times \Delta_1 &\to \Delta_1\\
            (a,b) &\mapsto a+b-ab.
          \end{aligned}
        \]
        Since $I_{\sA}^*$ preserves limits, it follows that the interval
        $(I_{\sA}^*(\Delta_1),I_{\sA}^*(e_0),I_{\sA}^*(e_1))$, equipped with the
        image by $I_{\sA}^*$ of the above binary operation, is multiplicative.
      \end{ex}
      \begin{lemma}\label{lemma:multlocasph}
        Let $\sM$ be a category with finite products, $\cW$ a weakly saturated
        class of maps of $\sM$, $(\bI,i_0,i_1)$ a interval in $\sM$ such that
        $\bI \to e_{\sM}$ is universally in $\cW$ and
        $(\bL,\lambda_0,\lambda_1,\Lambda)$ a multiplicative interval in $\sM$.
        If there exists a morphism of intervals $(\bI,i_0,i_1) \to
        (\bL,\lambda_0,\lambda_1)$, then $\bL \to e_{\sM}$ is universally in
        $\cW$.
      \end{lemma}
      \begin{proof}
        See \cite[Lemme 1.4.10]{maltsiniotis2005theorie}.
      \end{proof}
      
      For the next definition, recall that a $(2,1)$\nbd{}category is a
      $2$\nbd{}category such that every $2$\nbd{}morphism is invertible (in
      other words, a $\Grpd$\nbd{}enriched category). Limits and colimits in a
      $(2,1)$\nbd{}category are the $\Grpd$\nbd{}enriched ones.
      \begin{defn}
        Let $\underline{\sM}$ be a $(2,1)$\nbd{}category with finite products
        (the terminal object is denoted by $e_{\sM}$) and an initial object
        $\varnothing$. An interval $(\bI,i_0,i_1)$ (of the underlying category
        of) $\underline{\sM}$ is said to be \emph{strongly separating} if for
        every $2$\nbd{}square of $\underline{\sM}$
        \[
          \begin{tikzcd}
            X \ar[r] \ar[d] & e_{\sM} \ar[d]\ar[d,"i_0"] \\
            e_{\sM} \ar[r,"i_1"'] & \bI,
            \ar[from=2-1,to=1-2,Rightarrow,"\simeq",shorten <=0.5em, shorten
            >=0.5em]
          \end{tikzcd}
        \]
        we necessarily have $X=\varnothing$ and the $2$\nbd{}morphism is the
        identity.
      \end{defn}
      \begin{ex}
        Consider $\Cat$ as a $(2,1)$\nbd{}category, where the $2$\nbd{}morphisms
        are the natural isomorphisms between functors. Then $(\Delta_1,e_0,e_1)$
        is strongly separating.
      \end{ex}
      \begin{paragr}Let $\sA$ be a small category. The category $\pgrpd{\sA}$
        has a canonical structure of a $2$\nbd{}category where the
        $2$\nbd{}morphisms are the strict natural $2$\nbd{}transformations. That
        is, given two parallel morphisms $\varphi,\psi \colon X \to Y$ of
        $\pgrpd{\sA}$, a $2$\nbd{}morphism $\alpha \colon \varphi \Rightarrow
        \psi$ consists of a family of natural transformations
        \[
          \begin{tikzcd}
            X(a) \ar[r,bend left,"\varphi_a",""{name=A,below}] \ar[r,bend
            right,"\psi_a"',""{name=B,above}]&Y(a),
            \ar[from=A,to=B,Rightarrow,"\alpha_a"]
          \end{tikzcd}
        \]
        such that, for every $f \colon a \to a'$ in $\sA$, the following
        naturality condition is satisfied
        \[
          \begin{tikzcd}
            X(a') \ar[r,"X(f)"]&X(a) \ar[r,bend
            left,"\varphi_a",""{name=A,below}] \ar[r,bend
            right,"\psi_a"',""{name=B,above}]&Y(a)
            \ar[from=A,to=B,Rightarrow,"\alpha_a"]
          \end{tikzcd}
          =
          \begin{tikzcd}
            X(a') \ar[r,bend left,"\varphi_{a'}",""{name=A,below}] \ar[r,bend
            right,"\psi_{a'}"',""{name=B,above}]&Y(a')\ar[r,"Y(f)"] & Y(a).
            \ar[from=A,to=B,Rightarrow,"\alpha_{a'}"]
          \end{tikzcd}
        \]
        Notice that since $X$ and $Y$ take values in groupoids, every $\alpha_a$
        is invertible, and it follows that every $2$\nbd{}morphism of
        $\pgrpd{\sA}$ is also invertible. Hence, $\pgrpd{\sA}$ is a
        $(2,1)$\nbd{}category.
      \end{paragr}
      \begin{lemma}\label{lemma:stronglyseparating}
        Let $i \colon \sA \to \Cat$ be a functor, with $\sA$ a small category,
        such that for every object $a$ of $\sA$, the category $i(a)$ is not
        empty. Then, the interval $(I^*(\Delta_1),I^*(e_0),I^*(e_1))$ of
        $\pgrpd{\sA}$ is strongly separating.
      \end{lemma}
      \begin{proof}
        It is obvious that $I^* \colon \Cat \to \pgrpd{\sA}$ can be extended to
        a $(2,1)$\nbd{}functor, and so is its left adjoint $I_!$ which was
        defined in the proof of \cref{lemma:I*adjoint}. We obtain this way a
        $(2,1)$\nbd{}adjunction. The fact that $(\Delta_1,e_0,e_1)$ is strongly
        separating in $\Cat$ can be expressed as the fact that the commutative
        square
        \[
          \begin{tikzcd}
            \varnothing \ar[r] \ar[d] & e \ar[d]\ar[d,"e_0"] \\
            e \ar[r,"e_1"'] & \Delta_1,
          \end{tikzcd}
        \]
        is a $(2,1)$\nbd{}comma square.\footnote{The definition is the same as
          the usual notion of comma square in a $2$\nbd{}category, except every
          $2$\nbd{}morphism involved is invertible.} Since $I^*$ preserves
        $\Grpd$\nbd{}enriched limits, the following square
        \[
          \begin{tikzcd}
            \varnothing \ar[r] \ar[d] & \ast \ar[d]\ar[d,"I^*(e_0)"] \\
            \ast \ar[r,"I^*(e_1)"'] & I^*(\Delta_1)
          \end{tikzcd}
        \]
        is also a $(2,1)$\nbd{}comma square (we used that
        $I^*(\varnothing)\simeq \varnothing$, which follows from the hypothesis
        on the non-emptiness of the categories $i(a)$), which means exactly that
        the interval $(I^*(\Delta_1),I^*(e_0),I^*(e_1))$ is strongly separating.
        Details are left to the reader.
      \end{proof}
      \begin{lemma}\label{lemma:lawvereweakinitial}
        Let $\sA$ be a small category. For every strongly separating interval
        $(\bI,i_0,i_1)$ of $\pgrpd{\sA}$, there exists a morphism of intervals
        $(\bI,i_0,i_1) \to (I_{\sA}^*(\Delta_1),I_{\sA}^*(e_0),I_{\sA}^*(e_1))$
        (non-necessarily unique).
      \end{lemma}
      \begin{proof}
        By adjunction, we need to find a morphism of $\Cat$, $u \colon
        I_{\sA}(\bI) \to \Delta_1$, such that the following diagram is
        commutative
        \begin{equation}\label{diag:sieve}
          \begin{tikzcd}
            \sA \ar[d] \ar[r] & e \ar[d,"e_0"] \\
            I_{\sA}(\bI) \ar[r,"u"] & \Delta_1 \\
            \sA \ar[u] \ar[r] & e, \ar[u,"e_1"']
          \end{tikzcd}
        \end{equation}
        where the map $\sA \to I_{\sA}(\bI)$ at the top is defined as $a \mapsto
        (a,a\to \ast \overset{i_0}{\to} \bI)$, and the other one similarly with
        $i_1$ instead of $i_0$.

        Let $(a,a \overset{p}{\to} \bI)$ be an object of $I_{\sA}(\bI)$.
        \begin{itemize}[label=-]
        \item If $p$ is such that there exists a natural isomorphism
          \[
            \begin{tikzcd}
              a \ar[r] \ar[rd,"p"',""{name=A,above}] &\ast \ar[d,"i_0"] \\
              & \bI, \ar[from=A,to=1-2,Rightarrow,"\simeq"]
            \end{tikzcd}
          \]
          then we define $u(a,p)=0$,
        \item else we define $u(a,p)=1$.
        \end{itemize}
        Given a morphism $(a',p') \to (a,p)$ of $I_{\sA}(\bI)$, notice that if
        $u(a,p)=0$, then $u(a',p')=0$ too. (In other words, the objects sent to
        $0$ form a \emph{sieve}). This allows for a unique possible way of
        defining $u$ on arrows.

        The upper square of \eqref{diag:sieve} is commutative by definition. For
        the lower square, we need to prove that
        \[
          u(a,a\to \ast \overset{i_1}{\to} \bI)=1.
        \]
        Suppose that it is not the case: this would mean that there exists a
        $2$\nbd{}square
        \[
          \begin{tikzcd}
            a \ar[r] \ar[d] & \ast \ar[d]\ar[d,"i_0"] \\
            \ast \ar[r,"i_1"'] & \bI,
            \ar[from=2-1,to=1-2,Rightarrow,"\simeq",shorten <=0.5em, shorten
            >=0.5em]
          \end{tikzcd}
        \]
        which is forbidden since $(\bI,i_0,i_1)$ is strongly separating (and the
        initial presheaf $\varnothing$ is never representable).
      \end{proof}
    
      \begin{prop}
        Let $\sA$ be a small category. The following are equivalent:
        \begin{enumerate}[label={\upshape(\alph*)}]
        \item $\sA$ is groupoidal local test,
        \item $I_{\sA}^*(\Delta_1)$ is locally aspherical,
        \item there exists a strongly separating interval $(\bI,i_0,i_1)$ in
          $\pgrpd{\sA}$, such that $\bI$ is locally aspherical,
        \item there exists a groupoidal locally aspherical functor $i \colon \sA
          \to \Cat$.
        \end{enumerate}
      \end{prop}
      \begin{proof}
        By \cref{rem:grpdlocaltestfun}, $\sA$ is groupoidal local test if and
        only if $\sA \to \Cat, a \mapsto \sA/a$ is a groupoidal locally
        aspherical functor. Hence, the implication $(a) \Rightarrow (d)$ is
        trivial and the equivalence $(a)\Leftrightarrow (b)$ follows from
        \cref{prop:criteriongrpdasphfun}. From \cref{lemma:stronglyseparating},
        we know that $(I_{\sA}^*(\Delta_1),I_{\sA}^*(e_0),I_{\sA}^*(e_1))$ is a
        strongly separating interval, hence the implication $(b) \Rightarrow
        (c)$. For the implication $(c) \Rightarrow (b)$, we know from
        \cref{ex:multinterval} that
        $(I_{\sA}^*(\Delta_1),I_{\sA}^*(e_0),I_{\sA}^*(e_1))$ is a
        multiplicative interval. Then, \cref{lemma:lawvereweakinitial} implies
        that there exists a morphism of intervals $(\bI,i_0,i_1) \to
        (I_{\sA}^*(\Delta_1),I_{\sA}^*(e_0),I_{\sA}^*(e_1))$. Since, by
        hypothesis, $\bI$ is locally aspherical, it follows from
        \cref{lemma:multlocasph} that $I_{\sA}^*(\Delta_1)$ is locally
        aspherical.

        So far, we have shown $(c) \Leftrightarrow (b) \Leftrightarrow (a)
        \Rightarrow (d)$. Let us conclude with the implication $(d) \Rightarrow
        (c)$. If $i \colon \sA \to \Cat$ is a groupoidal locally aspherical
        functor, then, by definition, $I^*(\Delta_1)$ is locally aspherical.
        Besides, each category $i(a)$ has a terminal object and in particular is
        not empty, hence \cref{lemma:stronglyseparating} applies and
        $(I^*(\Delta_1),I^*(e_0),I^*(e_1))$ is a strongly separating interval.
      \end{proof}
      
      \section{Groupoidal strict test categories}
      \emph{We fix once and for all in this section a basic localizer $\cW$ of
        $\Cat$.}
      \begin{paragr}
        Recall that a small category $\sA$ is \emph{totally aspherical} if
        \begin{rome}
        \item $\sA$ is aspherical,
        \item the diagonal functor
          \[
            \delta \colon \sA \to \sA \times \sA
          \]
          is aspherical.
        \end{rome}
      \end{paragr}
      \begin{ex}\label{ex:totasphcat}
        A small category that has finite products (including the empty product)
        is totally aspherical \cite[Exemple 1.6.4]{maltsiniotis2005theorie}.
      \end{ex}
      \begin{defn}
        A small category $\sA$ is \emph{$\cW$\nbd{}groupoidal strict test}, or
        simply \emph{groupoidal strict test}, if the following conditions are
        satisfied
        \begin{enumerate}[label={\upshape(\alph*)}]
        \item $\sA$ is totally aspherical,
        \item $\sA$ is groupoidal test.
        \end{enumerate}
      \end{defn}
        
      In the following proposition, it is important to understand that
      ``finite'' includes ``empty''.
      \begin{prop}\label{prop:eqcondtotasphgrpd}
        Let $\sA$ be a small category. The following are equivalent:
        \begin{enumerate}[label={\upshape(\alph*)}]
        \item $\sA$ is totally aspherical,
        \item the functor $I_{\sA} \colon \pgrpd{\sA} \to \Cat$ preserves finite
          products up to weak equivalence, i.e.\ for every finite family
          $(X_i)_{i \in I}$ of objects of $\pgrpd{\sA}$, the canonical morphism
          \[
            I_{\sA}(\prod_{i \in I}X_i) \to \prod_{i \in I}I_{\sA}(X_i)
          \]
          is a weak equivalence,
        \item the class of aspherical objects of $\pgrpd{\sA}$ is stable by
          finite products, i.e.\ if $(X_i)_{i \in I}$ is a finite family of
          aspherical objects of $\pgrpd{\sA}$, then
          \[
            \prod_{i \in I}X_i
          \]
          is also aspherical,

        \item for every finite family of $(a_i)_{i \in I}$ of objects of $\sA$,
          seen as representable presheaves (and as objects of $\pgrpd{\sA}$ via
          the canonical inclusion $\psh{\sA} \hookrightarrow \pgrpd{\sA}$), the
          product in $\pgrpd{\sA}$
          \[
            \prod_{i \in I}a_i
          \]
          is aspherical.
        \end{enumerate}
      \end{prop}
      \begin{proof}
        Let us begin with $(a) \Rightarrow (b)$. For the empty product, this is
        simply saying that $I_{\sA}(\ast)\simeq \sA \to e$ is a weak
        equivalence, which is the case because a totally aspherical category is
        in particular aspherical. Now let $X$ and $Y$ be two objects of
        $\pgrpd{\sA}$, and notice that the following square
        \[
          \begin{tikzcd}
            I_{\sA}(X\times Y) \ar[d,"\zeta_{X\times Y}"'] \ar[r] &
            I_{\sA}(X)\times I_{\sA}(Y)
            \ar[d,"\zeta_X\times\zeta_Y"] \\
            \sA \ar[r,"\delta"'] & \sA \times \sA
          \end{tikzcd}
        \]
        is commutative and a pullback square. Since a product of Grothendieck
        fibrations is again a Grothendieck fibration, the right vertical arrow
        is a Grothendieck fibration and by \cref{pullbackasph} we deduce that
        the top horizontal arrow is aspherical. The general case follows from an
        immediate induction and the fact that weak equivalences in $\Cat$ are
        stable by finite products \cite[Proposition
        2.1.3]{maltsiniotis2005theorie}.

        The implication $(b) \Rightarrow (c)$ is immediate because a finite
        product of aspherical categories is aspherical. The implication $(c)
        \Rightarrow (d)$ is trivial.

        Finally, for the implication $(d) \Rightarrow (a)$, notice first that
        condition $(d)$ applied to the empty product gives that the terminal
        $\ast$ object of $\pgrpd{\sA}$ is aspherical, which means exactly that
        $\sA$ is aspherical as usual. Now let $a$ and $b$ the two objects of
        $\sA$, seen as representable presheaves and thus as objects of
        $\pgrpd{\sA}$. It is straightforward to check that
        \[
          I_{\sA}(a\times b) \simeq A/(a,b),
        \]
        where on the right hand side, $(a,b)$ has to be understood as an object
        of $\sA\times\sA$ and the slice is relative the diagonal functor $\delta
        \colon \sA \to \sA\times\sA$. This slice category being aspherical for
        every $(a,b) \in \sA\times\sA$ means exactly that $\delta$ is
        aspherical.
      \end{proof}
      Now, the crucial result is the following.
      \begin{lemma}\label{lemma:aspheqlocasph}
        Let $\sA$ be a totally aspherical category and $X$ an object of
        $\pgrpd{\sA}$. The following conditions are equivalent:
        \begin{enumerate}[label={\upshape(\alph*)}]
        \item $X$ is aspherical,
        \item $X$ is locally aspherical.
        \end{enumerate}
      \end{lemma}
      \begin{proof}
        A totally aspherical category being in particular aspherical, the
        implication $(b) \Rightarrow (a)$ has already been proved in
        \cref{cor:locasphisasph}. For the other implication, let $a$ be an
        object of $\sA$, which we see as a representable presheaf (and then as
        an object of $\pgrpd{\sA}$ via the canonical inclusion $\psh{\sA}
        \hookrightarrow \pgrpd{\sA}$). Thanks to \cref{prop:eqcondtotasphgrpd}
        and because representable presheaves are always aspherical, we know that
        \[
          a\times X
        \]
        is an aspherical object of $\pgrpd{\sA}$, which proves that $X$ is
        locally aspherical by \cref{lemma:characgrpdlocasph}. \qedhere
      \end{proof}
      The following results are straightforward consequences of the previous
      lemma.
      \begin{prop}
        Let $\sA$ be a small category and $i \colon \sA \to \Cat$ a functor. If
        $\sA$ is totally aspherical, then the following conditions are
        equivalent:
        \begin{enumerate}[label={\upshape(\alph*)}]
        \item $i$ is a groupoidal locally aspherical functor,
        \item $i$ is a groupoidal aspherical functor.
        \end{enumerate}
      \end{prop}
      \begin{prop}
        Let $\sA$ be a small category. If $\sA$ is totally aspherical, then the
        following are equivalent:
        \begin{enumerate}[label={\upshape(\alph*)}]
        \item $\sA$ is groupoidal strict test,
        \item $\sA$ is groupoidal test,
        \item $\sA$ is groupoidal weak test,
        \item $I_{\sA}(\Delta_1)$ is aspherical,
        \item there exists a strongly separating interval $(\bL,i_0,i_1)$ in
          $\pgrpd{\sA}$ such that $\bL$ is aspherical,
        \item there exists a groupoidal aspherical functor $i \colon \sA \to
          \Cat$.
        \end{enumerate}
      \end{prop}
      \section{Test categories vs.\ Groupoidal test categories}
      \begin{paragr}
        The comparison of the theory of groupoidal test categories and test
        categories relies on the following trivial but essential observation. If
        $\sD$ is a (small) category with no non-trivial isomorphisms, then for
        any (small) category $\sC$, the groupoid $\Homi_{\Cat}(\sC,\sD)$ doesn't
        have any non-trivial morphisms. In other words, $\Homi_{\Cat}(\sC,\sD)$
        is a set and we have
        \[
          \Homi_{\Cat}(\sC,\sD)=\Hom_{\Cat}(\sC,\sD).
        \]
        In particular, let $i \colon \sA \to \Cat$ be a functor, with $\sA$ a
        small category. For any category $\sD$ with no non-trivial isomorphisms,
        we have
        \[
          I^*(\sD)=i^*(\sD).
        \]
        (Remember that we consider $\psh{\sA}$ as a full subcategory of
        $\pgrpd{\sA}$.)
      \end{paragr}
      We then immediately have the following result.
      \begin{prop}\label{prop:complocasphfun}
        Let $i \colon \sA \to \Cat$ a functor such that for every object $a$ of
        $\sA$, the category $i(a)$ has a terminal object. The following are
        equivalent:
        \begin{enumerate}[label={\upshape(\alph*)}]
        \item $i$ is a groupoidal locally aspherical functor,
        \item $i$ is a locally aspherical functor.
        \end{enumerate}
      \end{prop}
   
      \begin{proof}
        Thanks to \cref{prop:criteriongrpdasphfun}, condition $(a)$ is
        equivalent to $I^*(\Delta_1)$ being locally aspherical. And condition
        $(b)$ means that $i^*(\Delta_1)$ is locally aspherical (see
        \cref{def:asphfun}). Since $\Delta_1$ has no non-trivial isomorphism, we
        have
        \[
          I^*(\Delta_1)=i^*(\Delta_1).
        \]
        To conclude, let us prove that a $\Set$\nbd{}valued presheaf is locally
        aspherical as an object of $\psh{\sA}$ (\cref{defn:wepsh}) if and only
        if it is locally aspherical as an object of $\pgrpd{\sA}$
        (\cref{def:grpdlocallyasph}). First notice that the canonical inclusion
        $\psh{\sA} \hookrightarrow \pgrpd{\sA}$ preserves and reflects limits,
        and preserve and reflects weak equivalences. Hence, given an object $X$
        of $\psh{\sA}$, if $X \to \ast$ is a universal weak equivalence of
        $\pgrpd{\sA}$, then it is also a universal weak equivalence of
        $\psh{\sA}$. The latter is the definition of locally aspherical object
        of $\psh{\sA}$ and the former is a characterisation of locally
        aspherical object of $\pgrpd{\sA}$ (\cref{lemma:characgrpdlocasph}(d)).
        This proves the ``if'' part. Conversely, if $X \to \ast$ is a universal
        weak equivalence of $\psh{\sA}$, then for every object $a$ of $\sA$,
        seen as a representable presheaf, the canonical projection $a\times X\to
        a$ is a weak equivalence of $\psh{\sA}$. Since $a$ is an aspherical
        object of $\psh{\sA}$ (because $i_{\sA}(a)=\sA/a$ has a terminal
        object), it follows that $a\times X$ is an aspherical object of
        $\psh{\sA}$. By \cref{rem:compatibleasph}, we deduce that it is also an
        aspherical object of $\pgrpd{\sA}$, which proves that $X$ is locally
        aspherical as an object of $\pgrpd{\sA}$ by
        \cref{lemma:characgrpdlocasph}(c).
      \end{proof}
      From this, we deduce our first comparison theorem.
      \begin{theorem}\label{firstcompthm}
        Let $\sA$ be a small category. We have the following equivalences:
        \begin{enumerate}[label={\upshape(\alph*)}]
        \item $\sA$ is groupoidal local test $\Leftrightarrow$ $\sA$ is local
          test,
        \item $\sA$ is groupoidal test $\Leftrightarrow$ $\sA$ is test,
        \item $\sA$ is groupoidal strict test $\Leftrightarrow$ $\sA$ is strict
          test.
        \end{enumerate}
      \end{theorem}
      \begin{proof}
        By \cref{rem:grpdlocaltestfun} (resp.\ \cref{rem:localtestfun}), $\sA$
        is a groupoidal local test category (resp.\ local test category) if and
        only if $\sA \to \Cat, a \mapsto \sA/a$ is a locally groupoidal
        aspherical functor (resp.\ locally aspherical functor). The equivalence
        (a) follows then from \cref{prop:complocasphfun}.

        By \cref{prop:localtestandasph} (resp.\ \cref{prop:crittestcat}), we
        know that $\sA$ is groupoidal test (resp.\ test) if and only if it is
        groupoidal locally test (resp.\ locally test) and aspherical. Hence, the
        equivalence (b) follows trivially from (a).

        Finally, $\sA$ is groupoidal strict test (resp.\ strict test) if it is
        groupoidal test (resp.\ test) and totally aspherical. Hence, the
        equivalence (c) follows trivially from (b).
      \end{proof}
      \begin{cor}\label{cor:inclusionDK}
        If $\sA$ is a test category (or equivalently a groupoidal test
        category), the canonical inclusion functor $\psh{\sA} \hookrightarrow
        \pgrpd{\sA}$ induces an equivalence at the level of homotopy categories
        \[
          \Ho(\psh{\sA}) \simeq \Ho(\pgrpd{\sA}).
        \]
      \end{cor}
      \begin{proof}
        Consider the commutative triangle
        \[
          \begin{tikzcd}
            \psh{\sA}\ar[dr,"i_{\sA}"'] \ar[r,hook] &\pgrpd{\sA}
            \ar[d,"I_{\sA}"]
            \\
            &\Cat.
          \end{tikzcd}
        \]
        If $\sA$ is a test category (or equivalently a groupoidal test
        category), then both vertical arrow of the previous triangle induce
        equivalences at the level of homotopy categories. The result follows
        then by a 2-out-of-3 property for equivalences of categories.
      \end{proof}
      \begin{rem}\label{rem:DKeq}
        In fact, the proof of the previous corollary is straightforwardly
        generalized to deduce that if $\sA$ is a test category, then $\psh{\sA}
        \hookrightarrow \pgrpd{\sA}$ induces a Dwyer--Kan equivalence
        \cite{barwick2012characterization} $(\psh{\sA},\cW_{\psh{\sA}})
        \overset{\sim}{\rightarrow} (\pgrpd{\sA},\cW_{\pgrpd{\sA}})$, hence an
        equivalence of $(\infty,1)$\nbd{}categories.
      \end{rem}
      \begin{ex}
        It follows from \cref{firstcompthm}, that all the examples of (strict)
        test categories given in \cref{ex:testcategories} are also groupoidal
        (strict) test categories. In particular, $\Delta$ is a groupoidal strict
        test category, and we recover that the classical result that the
        category $\pgrpd{\Delta}$ models homotopy types (using
        \cref{cor:inclusionDK} for example). But this is also the case of
        $\Grpd$\nbd{}valued presheaves over: the cube category with or without
        connections, Joyal's $\Theta$ category, the dendroidal category, etc.
      \end{ex}
      Let us now compare groupoidal weak test categories with weak test
      categories. For this, first recall the following technical result.
      \begin{lemma}\label{lemma174}
        Let $u \colon \sA \to \sB$ be a morphism of $\Cat$, and $i \colon \sA
        \to \Cat$ and $j \colon \sB \to \Cat$ such that the triangle
        \[
          \begin{tikzcd}
            \sA \ar[r,"u"] \ar[dr,"i"']& \sB \ar[d,"j"] \\
            &\Cat
          \end{tikzcd}
        \]
        is commutative, and suppose that for every object $b$ of $\sB$, the
        category $j(b)$ has a terminal object.
        \begin{enumerate}[label={\upshape(\alph*)}]
        \item if $u$ is an aspherical morphism of $\Cat$, then $i \colon \sA \to
          \Cat$ is an aspherical functor if and only if $j \colon \sB \to \Cat$
          is,
        \item if $j \colon \sB \to \Cat$ is fully faithful and $i \colon \sA \to
          \Cat$ is an aspherical functor, then $u$ is an aspherical morphism of
          $\Cat$ and $j \colon \sB \to \Cat$ is an aspherical functor.
        \end{enumerate}
      \end{lemma}
      \begin{proof}
        See \cite[Lemma 1.7.4]{maltsiniotis2005theorie}.
      \end{proof}
      For $\Grpd$\nbd{}valued presheaves, we have the following partial
      generalization.
      \begin{lemma}\label{lemma174grpd}
        Let $u \colon \sA \to \sB$ be a morphism of $\Cat$, and $i \colon \sA
        \to \Cat$ and $j \colon \sB \to \Cat$ such that the triangle
        \[
          \begin{tikzcd}
            \sA \ar[r,"u"] \ar[dr,"i"']& \sB \ar[d,"j"] \\
            &\Cat
          \end{tikzcd}
        \]
        is commutative, and suppose that for every object $b$ of $\sB$, the
        category $j(b)$ has a terminal object.
        \begin{enumerate}[label={\upshape(\alph*)}]
        \item if $u$ is an aspherical morphism of $\Cat$, then $i \colon \sA \to
          \Cat$ is a groupoidal aspherical functor if and only if $j \colon \sB
          \to \Cat$ is,
        \item if $j \colon \sB \to \Cat$ is fully faithful, $i \colon \sA \to
          \Cat$ is a groupoidal aspherical functor, and for every object $b$ of
          $\sB$, the category $j(b)$ does not have any non-trivial isomorphism,
          then $u$ is an aspherical morphism of $\Cat$ and $j \colon \sB \to
          \Cat$ is a groupoidal aspherical functor.
        \end{enumerate}
      \end{lemma}
      \begin{proof}
        Notice first that the hypotheses imply that the category $i(a)$ has a
        terminal object for every object $a$ of $\sA$.
          
        Now, the given commutative triangle induces a commutative triangle
        \[
          \begin{tikzcd}
            \pgrpd{\sA}  & \pgrpd{\sB} \ar[l,"u^*"']\\
            & \Cat. \ar[lu,"I^*"] \ar[u,"J^*"']
          \end{tikzcd}
        \]
        If $u$ is aspherical, it follows from \cref{prop:asphmorgrpd} that for a
        small aspherical category $\sC$, $J^*(\sC)$ is aspherical if and only if
        $I^*(\sC)$ is aspherical. In particular, this proves $(a)$.

        For $(b)$, notice that with the hypotheses, we have, for every object
        $b$ of $\sB$,
        \[
          J^*(j(b))=j^*(j(b))\simeq b,
        \]
        where the first equality comes from the fact that $j(b)$ does not have
        any non-trivial isomorphism and the second from the fact that $j$ is
        fully faithful (note that on the right hand side of the second equality,
        we abusively wrote $b$ for the $\Set$\nbd{}valued presheaf represented
        by $b$). We then have
        \[
          u^*(b)\simeq u^*(J^*(j(b)))=I^*(j(b)).
        \]
        Since $i$ is a groupoidal aspherical functor and $j(b)$ is an aspherical
        category, we have that $I^*(j(b))$ is aspherical and so is $u^*(b)$.
        This means that the category $i_{\sA}(u^*(b))$ is aspherical and an
        immediate verification shows that we have a canonical isomorphism
        $i_{\sA}(u^*(b))\simeq A/b$, which proves by definition that $u$ is
        aspherical. Hence, we can apply (a) and the conclusion follows.
      \end{proof}
      \begin{rem}\label{rem:nontrivialiso}
        As the reader might notice, it is only the $(b)$ of \cref{lemma174} that
        does not generalize straightforwardly in \cref{lemma174grpd} and for
        which we added the hypothesis that for every $b$ in $\sB$, $j(b)$ does
        not have non-trivial isomorphisms. We do not know whether this
        hypothesis is necessary or not.
      \end{rem}
   
      We can now prove the following result.
      \begin{prop}\label{prop:asphfunvsgrpdasphfun}
        Let $i \colon \sA \to \Cat$ be a functor, with $\sA$ a small category,
        such that $i(a)$ has a terminal object for every object $a$ of $\sA$. We
        have the following implication:
        \begin{center}
          $i$ is an aspherical functor $\Rightarrow$ $i$ is a groupoidal
          aspherical functor.
        \end{center}
        If we suppose moreover that for every object $a$ of $\sA$, the category
        $i(a)$ does not have any non-trivial isomorphism, then we also have the
        converse implication:
        \begin{center}
          $i$ is a groupoidal aspherical functor $\Rightarrow$ $i$ is an
          aspherical functor.
        \end{center}
      \end{prop}
      \begin{proof}
        We begin by the first implication. Suppose that $i \colon \sA \to \Cat$
        is an aspherical functor and let $\sB$ be the smallest full subcategory
        of $\Cat$ such that:
        \begin{itemize}[label=-]
        \item every $i(a)$, for $a$ in $\sA$, is an object of $\sB$,
        \item $\sB$ is stable by finite products.
        \end{itemize}
        Since a finite product of categories with terminal object has a terminal
        object, it follows that every object of $\sB$ has a terminal object. By
        construction, we have a factorization,
        \[
          \begin{tikzcd}
            \sA \ar[r,"i_0"] \ar[dr,"i"'] & \sB \ar[d,"j",hook]\\
            &\Cat,
          \end{tikzcd}
        \]
        where $j$ is the canonical inclusion and $i_0 \colon \sA \to \sB$ is the
        functor $a \mapsto i(a)$. Since $i$ is an aspherical functor,
        \cref{lemma174}(b) implies that $i_0$ is aspherical and $j \colon \sB
        \hookrightarrow \Cat$ is an aspherical functor. Since $\sB$ is stable by
        finite products, it is totally aspherical (\cref{ex:totasphcat}), and
        thus $j$ is also a locally aspherical functor (this follows easily from
        \cref{lemma:aspheqlocasph} applied to $\Set$\nbd{}valued presheaves).
        Applying \cref{prop:complocasphfun}, we obtain that $j$ is a groupoidal
        locally aspherical functor, and then a groupoidal aspherical functor
        because $\sB$ is aspherical. Finally, using \cref{lemma174grpd}(a), we
        have that $i$ is a groupoidal aspherical functor.

        For the converse implication, let $i \colon \sA \to \Cat$ be a
        groupoidal aspherical functor such that each $i(a)$, for a in $\sA$,
        does not have any non-trivial isomorphism and let $\sB$ be the smallest
        full subcategory of $\Cat$ such that:
        \begin{itemize}[label=-]
        \item every $i(a)$, for $a$ in $\sA$, is an object of $\sB$,
        \item every object $\Delta_n$ of $\Delta$, for $n\geq 0$, is an object
          of $\sB$.
        \end{itemize}
        Notice that the objects of $\sB$ do not have any non-trivial
        isomorphism. By construction, we have a commutative diagram
        \[
          \begin{tikzcd}
            \sA \ar[r,"i_0"] \ar[dr,"i"'] & \sB \ar[d,"j",hook] & \Delta
            \ar[l,"k_0"'] \ar[dl,"k",hook]\\
            & \Cat,
          \end{tikzcd}
        \]
        where:
        \begin{itemize}[label=-]
        \item $i_0$ is the functor $a \mapsto i(a)$,
        \item $j$ is the full subcategory inclusion,
        \item $k$ is the canonical fully faithful inclusion of $\Delta$ in
          $\Cat$,
        \item $k_0$ is the full subcategory inclusion.
        \end{itemize}
        By hypothesis, $i$ is a groupoidal aspherical functor and since $j$ is
        fully faithful and the objects of $\sB$ do not have any non-trivial
        isomorphism, it follows from \cref{lemma174grpd}(b) that $i_0$ is an
        aspherical morphism of $\Cat$. Now, since $k \colon \Delta \to \Cat$ is
        an aspherical functor \cite[Exemple 1.17.18]{maltsiniotis2005theorie}
        (the functor $k^* \colon \Cat \to \psh{\Delta}$ is nothing but the nerve
        functor), and since $j$ is fully faithful, it follows from
        \cref{lemma174}(b) that $j$ is also an aspherical functor. Using that
        $i_0$ is an aspherical morphism of $\Cat$, we deduce from an application
        \cref{lemma174}(a) that $i$ is an aspherical functor.
      \end{proof}
      From this we deduce the following comparison theorem.
      \begin{theorem}
        Let $\sA$ be a small category. We have the following implication:
        \begin{center}
          $\sA$ is a weak test category $\Rightarrow$ $\sA$ is a groupoidal weak
          test category.
        \end{center}
        If we suppose moreover that $\sA$ does not have any non-trivial
        isomorphism, then we also have the converse implication:
        \begin{center}
          $\sA$ is a groupoidal weak test category $\Rightarrow$ $\sA$ is a weak
          test category.
        \end{center}
      \end{theorem}
      \begin{proof}
        By \cref{prop:critgrpdweaktest} (resp.\ \cref{rem:localtestfun}), $\sA$
        is groupoidal weak test (resp.\ weak test) if and only if $\sA \to \Cat,
        a \mapsto \sA/a$ is a groupoidal aspherical functor (resp. aspherical
        functor). Hence, the result follows immediately from
        \cref{prop:asphfunvsgrpdasphfun}.
      \end{proof}
      \begin{rem}
        In light of \cref{rem:nontrivialiso}, we do not know if the hypothesis
        that $\sA$ does not have any non-trivial isomorphism is necessary for
        the second implication of the previous theorem. If a counter-example
        exists, then it would necessary be a small category $\sA$ with
        non-trivial isomorphisms which is groupoidal weak test but not
        groupoidal test (or else \cref{firstcompthm} applies).
      \end{rem}
      \begin{cor}
        If $\sA$ is a weak test category, then the canonical inclusion
        $\psh{\sA} \to \pgrpd{\sA}$ induces an equivalence at the level of
        homotopy categories
        \[
          \Ho(\psh{\sA}) \simeq \Ho(\pgrpd{\sA}).
        \]
      \end{cor}
      \begin{proof}
        Similar to the proof of \ref{cor:inclusionDK}.
      \end{proof}
      \begin{rem}
        Same remark as \cref{rem:DKeq}.
      \end{rem}
      \begin{ex}
        It follows from \cref{firstcompthm} that all examples of weak test
        categories from \cref{ex:testcategories} are also groupoidal weak test
        categories. For example, it is the case of the category $\Delta'$ of
        finite non-empty ordinals and non-decreasing monomorphisms. In
        particular, the category $\pgrpd{(\Delta')}$ models homotopy types.
      \end{ex}
      \begin{paragr}
        Finally, let us end this section with a quick word on the comparison of
        pseudo-test categories and groupoidal pseudo-test categories. Although,
        the following implication
        \begin{equation}\label{conjpstest}
          \text{ pseudo-test category } \Rightarrow \text{ groupoidal pseudo-test category }
        \end{equation}
        seems reasonable to expect, it remains an open question for the author
        of these notes. As for the converse implication, it is not true in
        general. More precisely, the example below shows that there exists a
        basic localizer $\cW$ (which is not $\cW_{\infty}$ !) such that the
        class of $\cW$\nbd{}groupoidal pseudo-test categories strictly contains
        the class of $\cW$\nbd{}pseudo-test categories. (Hence, for this
        \emph{particular} basic localizer, the implication \eqref{conjpstest} is
        true.) The question remains open for an arbitrary basic localizer, in
        particular for~$\cW_{\infty}$.
      \end{paragr}
      \begin{ex}
        Consider the functor $\pi_1 \colon \Cat \to \Grpd$, left adjoint of the
        canonical inclusion functor $\iota \colon \Grpd \to \Cat$, and let
        $\cW_{1}$ be the class of morphisms $f$ of $\Cat$ such that $\pi_{1}(f)$
        is an equivalence of groupoids. We leave it as an exercise to the reader
        to show that $\cW_{1}$ is a basic localizer of $\Cat$. Now, since
        $\iota$ is fully faithful, the co--unit of the adjunction $\pi_1 \dashv
        \iota$ is an isomorphism and it follows then from
        \cref{lemma:crithomotopicalequiv} that this adjunction induces a
        homotopical equivalence between $(\Cat,\cW_{1})$ and
        $(\Grpd,\cW_{\mathrm{EqGrpd}})$, where $\cW_{\mathrm{EqGrpd}}$ is the
        class of equivalences of groupoids. (This proves in particular that
        $\Cat$ models homotopy $1$\nbd{}types.) It follows that a small category
        $\sA$ is $\cW_{1}$\nbd{}groupoidal pseudo-test if and only if it is
        $\cW_1$\nbd{}aspherical and the functor
        \[
          \begin{aligned}
            \pgrpd{\sA} &\to \Grpd\\
            X &\mapsto \pi_1(i_A(X))
          \end{aligned}
        \]
        induces an equivalence at the level of homotopy categories.

        Now, let $\sA=e$ be the terminal category. Then, the previous functor is
        nothing but the identity functor of $\Grpd$ and it follows trivially
        that $e$ is $\cW_1$\nbd{}groupoidal pseudo-test. On the other hand, $e$
        is not $\cW_{1}$\nbd{}pseudo-test. Indeed, $i_e \colon \Set \to \Cat$ is
        nothing but the canonical inclusion functor and so $i_e^{-1}(\cW_{1})$
        is the class of isomorphisms of $\Set$. If $e$ was $\cW_{1}$
        pseudo-test, this would imply that $(\Set,\mathrm{Iso}) \hookrightarrow
        (\Grpd,\cW_1)$ induces an equivalence at the level homotopy categories
        (and that $\Set$ models homotopy $1$\nbd{}types), which is easily seen
        to be false.
      \end{ex}
      \begin{rem}
        Note that the previous example also shows that the class of
        $\cW_1$\nbd{}groupoidal weak test categories is strictly bigger than the
        class of $\cW_1$\nbd{}groupoidal pseudo-test categories. Indeed, if the
        terminal category $e$ were a $\cW_1$\nbd{}groupoidal weak test category,
        then it would be a $\cW_1$\nbd{}weak test category (since it does not
        have any non-trivial isomorphisms), and in particular a
        $\cW_1$\nbd{}pseudo-test category.
      \end{rem}
      \section{Weak equivalences via the nerve}\label{sec:wenerve}
      The section of this paper is devoted to giving an equivalent definition of
      weak equivalences of $\Grpd$\nbd{}valued presheaves in terms of nerve
      functors. In particular, in \cref{ex:simplicialwe} below, we recover the
      usual definition of weak equivalences on $\pgrpd{\Delta}$ used in the
      literature \cite[section 8]{crans1995quillen}, \cite{joyal1996homotopy}.

      \begin{paragr}\label{paragr:GrC}
        Let $\sA$ be a small category. For a functor $X \colon \sA^{\op} \to
        \Cat$, we denote by $\int_{\sA}X$ the \emph{Grothendieck construction of
          $X$}. This means that $\int_{\sA}X$ is the category such that:
        \begin{itemize}[label=-]
        \item objects are pairs $(a,x)$ where $a$ is an object of $\sA$ and $x$
          is an object of $X(a)$,
        \item a morphism $(a,x) \to (a',x')$ is a pair $(f,k)$ where $f \colon a
          \to a'$ is a morphism of $\sA$ and $k \colon x \to X(f)(x')$ is a
          morphism of $X(a)$.
        \end{itemize}
        (For details, we refer to \cite[2.2.6]{maltsiniotis2005theorie}, where
        the notation $\nabla_{\sA}$ for $\int_{\sA}$ is used). This construction
        is functorial and provides a functor
        \[
          \int_{\sA} \colon \pcat{\sA} \to \Cat,
        \]
        where we write $\pcat{\sA}$ for the category of functors $\sA^{\op} \to
        \Cat$ and natural transformations between them.
            
        We have canonical inclusions $\psh{\sA} \hookrightarrow \pgrpd{\sA}
        \hookrightarrow \pcat{\sA}$, and, as already observed, for $X$ an object
        of $\psh{\sA}$ (resp. $\pgrpd{\sA}$), we have $\int_{\sA}X=i_{\sA}(X)$
        (resp.\ $\int_{\sA}X=I_{\sA}(X)$).
      \end{paragr}
      \begin{prop}\label{ThomColimit}\cite[Proposition
        2.3.1]{maltsiniotis2005theorie} Let $\cW$ be a basic localizer of $\Cat$
        and $\sA$ a small category. The functor $\int_{\sA} \colon \pcat{\sA}
        \to \Cat$ sends pointwise $\cW$\nbd{}equivalences to
        $\cW$\nbd{}equivalences.
      \end{prop}
      \emph{We now fix once and for all a basic localizer $\cW$ of $\Cat$.}
      \begin{paragr}
        Let $\sA$ and $\sB$ be two small categories and consider the presheaf
        category $\psh{\sA\times\sB}$. Using the identification
        $\psh{\sA\times\sB}\simeq \underline{\mathrm{Hom}}(\sA^{\op},\psh{\sB})$
        and the post-composition by the functor $i_{\sB} \colon \psh{\sB} \to
        \Cat$ defines a functor
        \[
          i_{\sB} \colon \psh{\sA\times\sB} \to \pcat{\sA}
        \]
        which we abusively denote by $i_{\sB}$ again. Then if we postcompose by
        $\int_{\sA}$, we obtain a functor
        \[
          \int_{\sA}i_{\sB} \colon \psh{\sA\times\sB} \to \Cat.
        \]
      \end{paragr}
      The proof of the following lemma is a straightforward verification, which
      we leave to the reader.
      \begin{lemma}\label{lemma:intmult}
        For every object $X$ of $\psh{\sA\times\sB}$, there is a canonical
        isomorphism
        \[
          i_{\sA\times\sB}(X) \simeq \int_{\sA}i_{\sB}(X),
        \]
        which is natural in $X$.
      \end{lemma}
      \begin{rem}
        Remember that the Grothendieck construction of a functor with values in
        $\Cat$ is weakly equivalent to its homotopy colimit (with respect to any
        basic localizer on $\Cat$) \cite[Théorème
        3.1.7]{maltsiniotis2005theorie}.\footnote{When the basic localizer is
          $\cW_{\infty}$, this is a result of Thomason
          \cite{thomason1979homotopy}.} Since the functor $i_{\sB}$ is just the
        restriction of the Grothendieck construction to Set-valued presheaves,
        the previous lemma can be simply restated by saying that the homotopy
        colimit of a functor of two variables is computed by successively taking
        the homotopy colimit relative to each variable.
      \end{rem}
      \begin{paragr}
        Let $\sA$ and $\sB$ be small categories and let $i \colon \sB \to \Cat$
        be a functor. Recall that we denote $i^* \colon \Cat \to \psh{\sB}$ the
        functor $\sC \mapsto \Hom_{\Cat}(i(-),\sC)$. By considering $\Grpd$ as a
        subcategory of $\Cat$, the functor $i^*$ induces by post-composition a
        functor
        \[
          i^* \colon \pgrpd{\sA} \to \psh{\sA \times \sB},
        \]
        which we abusively denote by $i^*$ as
        well. 
      \end{paragr}
      \begin{prop}Let $\sA$ and $\sB$ be small categories and $i \colon \sB \to
        \Cat$ a functor such that for every $b$ in $\sB$, the category $i(b)$
        has a terminal object. Then, there exists a natural transformation
        \[
          \begin{tikzcd}[column sep=small]
            \pgrpd{\sA} \ar[rr,"i^*"]\ar[dr,"I_{\sA}"',""{name=A,right}] &&
            \psh{\sA\times\sB}
            \ar[dl,"i_{\sA \times \sB}"] \\
            &\Cat.& \ar[from=1-3,to=A,Rightarrow,shorten <=1em, shorten >=1em]
          \end{tikzcd}
        \]
        Moreover, if $i$ is an aspherical functor, then this natural
        transformation is a weak equivalence argument by argument.
      \end{prop}
      \begin{proof}
        For every $b$ in $\sB$, let $e_{b}$ be the terminal object of $i(b)$. By
        an analogous construction as the one in \cref{yogatermobj} in the case
        of $\Set$\nbd{}valued presheaves (see
        \cite[3.2.4]{cisinski2006prefaisceaux} for details), for every small
        category $\sC$ we define a morphism
        \begin{equation}\label{alphaterm}
          \begin{aligned}
            \alpha_{\sC} \colon i_{\sB}i^*(\sC) &\to \sC \\
            (b,p \colon i(b) \to \sC) &\mapsto p(e_b),
          \end{aligned}
        \end{equation}
        which is natural in $\sC$. For every $X$ in $\pgrpd{\sA}$ and $a$ in
        $\sA$, we obtain a map
        \[
          \alpha_{X(a)} \colon i_{\sB}i^*(X(a)) \to X(a),
        \]
        natural in $X$ and $a$, and by applying $\int_{\sA}$, we obtain a
        canonical map
        \[
          i_{\sA\times \sB}i^*(X)\simeq \int_{\sA}i_{\sB}i^*(X) \longrightarrow
          \int_{\sA}X=I_{\sA}(X),
        \]
        which is natural in $X$. Now, if $i$ is an aspherical functor, then by
        \cite[Proposition 1.7.6]{maltsiniotis2005theorie} (which is the analogue
        of our \cref{prop:critgrpdasphfun} for $\Set$\nbd{}valued presheaves),
        the map \eqref{alphaterm} is a weak equivalence. We conclude with
        \cref{ThomColimit}.
      \end{proof}
      \begin{cor}\label{cor:wenerve}
        Let $\sA$ and $\sB$ be small categories, and $i \colon \sB \to \Cat$ an
        aspherical functor (such that $i(b)$ has a terminal object for every
        object $b$ of $\sB$). Then $i^* \colon \pgrpd{\sA} \to
        \psh{\sA\times\sB}$ preserves and reflects weak equivalences, i.e.
        \[
          \cW_{\pgrpd{\sA}}=i^{*-1}(\cW_{\psh{\sA\times\sB}}).
        \]
      \end{cor}
      A particular case where the previous corollary applies is the following.
      \begin{cor}
        Let $\sA$ be a totally aspherical small category and $i \colon \sA \to
        \Cat$ aspherical functor (such that $i(a)$ has a terminal object for
        every object $a$ of $\sA$). Then a morphism of $\pgrpd{\sA}$ is a weak
        equivalence if and only if its image by $i^* \colon \pgrpd{\sA} \to
        \psh{\sA\times\sA}$ is a diagonal weak equivalence, i.e.\
        \[
          \cW_{\pgrpd{\sA}}=i^{*-1}(\delta^{*-1}(\cW_{\psh{\sA}})),
        \]
        where $\delta^* \colon \psh{\sA\times\sA} \to \psh{\sA}$ is the diagonal
        functor.
      \end{cor}
      \begin{proof}
        If $\sA$ is totally aspherical, then the diagonal functor $\delta \colon
        \sA \to \sA\times\sA$ is aspherical and so the induced functor $\delta^*
        \colon \psh{\sA\times\sA} \to \psh{\sA}$ preserves and reflects weak
        equivalences \cite[Proposition 1.2.9(d)]{maltsiniotis2005theorie}.
      \end{proof}
      \begin{ex}\label{ex:simplicialwe}
        Let $\sA=\sB=\Delta$ and $i \colon \Delta \to \Cat$ be the canonical
        inclusion, so that $i^* \colon \Cat \to \psh{\Delta}$ is nothing but the
        usual nerve functor. The previous corollary implies that the weak
        equivalences on simplicial groupoids $\pgrpd{\Delta}$ are exactly those
        morphisms that induce diagonal weak equivalences of bisimplicial sets.
      \end{ex}
      \section{Bonus result: groupoids internal to categories model homotopy
        types}\label{section:grpdcatmodel}
      \begin{paragr}\label{paragr:internalgrpd}
        Let $\cC$ be a category with finite limits (or even only pullbacks). A
        groupoid internal to $\cC$ consists of a pair of objects $(X_0,X_1)$ of
        $\cC$ equipped with
        \begin{itemize}[label=-]
        \item source and target maps $s,t \colon X_1 \to X_0$,
        \item a unit map $i \colon X_0 \to X_1$, such that $s\circ i = t\circ i
          = \id_{X_0}$,
        \item an inverse map $inv \colon X_1 \to X_1$, such that $s \circ
          \mathrm{inv} = t$ and $t \circ \mathrm{inv} = s $,
        \item a composition map $m \colon X_1\times_{X_0}X_1 \to X_1$, such that
          $s \circ m = s \circ \pi_1$ and $t \circ m = t \circ \pi_2$, where
          $X\times_{X_0}X_1$, $\pi_1$ and $\pi_2$ are defined as the following
          fibred product
          \[
            \begin{tikzcd}
              X_1\times_{X_0}X_1 \ar[r,"\pi_1"] \ar[d,"\pi_2"] & X_1 \ar[d,"t"] \\
              X_1 \ar[r,"s"']& X_0, \ar[from=1-1,to=2-2,phantom,"\lrcorner",
              very near start]
            \end{tikzcd}
          \]
        \end{itemize}
        all of which satisfy the usual axioms saying that: $m$ is associative,
        $i$ is the unit on the left and right of the composition, and $inv$
        gives the inverse on the left and right of the composition.

        We shall often abuse notation and refer to a groupoid internal to $\cC$
        as a pair $X=(X_0,X_1)$
        
        An internal morphism of groupoids $f \colon X \to X'$ consists of a pair
        of morphisms $(f_0 \colon X_0 \to X_0', f_1 \colon X_1 \to X_1')$ of
        $\cC$ which commute with source, target, inverse, unit and composition
        in the obvious way. Internal groupoids and internal morphisms of
        groupoids form a category denoted by $\Grpd(\cC)$.

        If $\cC'$ is another category with pullbacks and $F \colon \cC \to \cC'$
        is a functor preserving pullbacks, then $F$ sends groupoids internal to
        $\cC$ to groupoids internal to $\cC'$, hence there is an induced functor
        \[
          \Grpd(F) \colon \Grpd(\cC) \to \Grpd(\cC').
        \]
        Furthermore, if $F$ admits a left adjoint $G \colon \cC' \to \cC$ that
        \emph{also preserves pullbacks}, then we have an induced adjunction
        $\Grpd(G)\dashv \Grpd(F)$. The unit and co-unit of this adjunction are
        obtained by applying those of the adjunction $G \dashv F$ levelwise.
        This means that if $(X_1,X_0)$ is an object of $\Grpd(\cC)$, then the
        co-unit of this adjunction is simply
        \[
          (\varepsilon_{X_1} \colon GF(X_1)\to X_1, \varepsilon_{X_0} \colon
          GF(X_0) \to X_0),
        \]
        where $\varepsilon$ is the co-unit of the adjunction $G \dashv F$, and
        similarly for the unit.
      \end{paragr}
      \begin{rem}
        Even if the left adjoint of $F$ does not preserve pullbacks, under mild
        conditions (e.g.\ $\cC$ and $\cC'$ are locally presentable and $F$ is
        accessible) the functor $\Grpd(F)$ still admits a left adjoint. In
        general, this left adjoint is hard to manipulate and is constructed
        abstractly. The point of the previous paragraph is that, \emph{when the
          left adjoint $G$ of $F$ preserves pullbacks}, then the left adjoint of
        $\Grpd(F)$ is simply $\Grpd(G)$, and furthermore the unit and co-unit
        have a particularly nice form. This will play an important role later in
        the proof of \cref{prop:grpdcatmodel} (cf.\ also
        \cref{rem:grpdcatmodel}).
      \end{rem}
      \begin{paragr}\label{paragr:nerve}
        Let $\cC$ be a category with pullbacks and $X=(X_0,X_1)$ an internal
        groupoid of $\cC$. The \emph{nerve of $X$} is the simplicial object
        $N_{\ast}(X) \colon \Delta^{\op} \to \cC$ defined by
        \[
          N_{n}(X)=\underbrace{X_{1}\times_{X_0}\cdots\times_{X_0}X_1}_{n\text{
              times}},
        \]
        where for $n\geq 2$, this means that $N_{n}(X)$ is the limit of the
        diagram
        \[
          \begin{tikzcd}[column sep=tiny]
            X_1 \ar[dr,"s"] && X_1 \ar[dl,"t"'] \ar[dr,"s"] & & \ar[dl,"t"']\cdots\ar[dr,"s"] && X_1 \ar[dl,"t"']\\
            &X_0&&X_0&&X_0,
          \end{tikzcd}
        \]
        where $X_1$ appears $n$ times, and, by convention, for $n=0,1$
        \[
          \begin{aligned}
            N_{0}(X)&=X_0,\\
            N_{1}(X)&=X_1.
          \end{aligned}
        \]
        The face and degeneracy maps are defined as follows:
        \begin{itemize}[label=-]
        \item for $0 < j < n$, $\partial_j \colon N_{n}(X) \to N_{n-1}(X)$ is
          induced by $m \colon X_1\times_{X_0}X_1 \to X_1$ acting on the
          $(j-1)$\nbd{}th and $j$\nbd{}th factors of $N_{n}(X)$,
        \item for $n\geq 1$, $j=0$ or $j=n$, $\partial_j \colon N_{n}(X) \to
          N_{n-1}(X)$ is the canonical projection that discards the $j$\nbd{}th
          factor of $N_{n}(X)$,
        \item for $0 \leq j \leq n$, $\epsilon_{j} \colon N_{n}(X) \to
          N_{n+1}(X)$ is induced by $i \colon X_0 \to X_1$ hitting the
          $j$\nbd{}th factor of $N_{n+1}(X)$.
        \end{itemize}
        Note that this construction does not use the inverse map $\mathrm{inv}
        \colon X_1 \to X_1$ and thus only depends on the underlying category
        (internal to $\cC$) of $X$.
      \end{paragr}
      \begin{lemma}\label{lemma:comtensprod}
        Let $\cC$, $\cC'$ be categories with pullbacks and $F \colon \cC \to
        \cC'$ a functor preserving pullbacks. The following diagram is
        commutative (up to an isomorphism of functors)
        \[
          \begin{tikzcd}[column sep=huge]
            \Grpd(\cC)\ar[d,"N_*"'] \ar[r,"\Grpd(F)"] & \Grpd(\cC') \ar[d,"N_*"]
            \\
            \Homu(\Delta^{\op},\cC) \ar[r,"{\Homu(\Delta^{\op},F)}"'] &
            \Homu(\Delta^{\op},\cC')
          \end{tikzcd}
        \]
      \end{lemma}
      \begin{proof}
        Straightforward verification left to the reader.
      \end{proof}
      \emph{We now fix once and for all a basic localizer $\cW$ on $\Cat$}.

      \begin{paragr}
        Let $\Grpd(\Cat)$ be the category of groupoids internal to $\Cat$.
        Recall that we denote by $\pcat{\Delta}$ the category of
        $\Cat$\nbd{}valued presheaves over $\Delta$. The construction from
        \cref{paragr:nerve} yields a functor $N_{\ast} \colon \Grpd(\Cat) \to
        \pcat{\Delta}$. If we post-compose by the Grothendieck construction, we
        obtain a functor from $\Grpd(\Cat)$ to $\Cat$:
        \[
          \Grpd(\Cat)\overset{N_{\ast}}{\longrightarrow} \pcat{\Delta}
          \overset{\int_{\Delta}}{\longrightarrow} \Cat.
        \]
      \end{paragr}
      \begin{defn}
        A morphism $f \colon X \to Y$ of $\Grpd(\Cat)$ is a \emph{weak
          equivalence} if
        \[
          \int_{\Delta}N_{\ast}(f) \colon \int_{\Delta}N_{\ast}(X) \to
          \int_{\Delta}N_{\ast}(Y)
        \]
        is in $\cW$.
      \end{defn}
      \begin{paragr}
        Let $\sA$ be a small category and consider the adjunction $i_{\sA}
        \colon \psh{\sA} \rightleftarrows \Cat \colon i_{\sA}^*$. The key
        observation is that the left adjoint $i_{\sA}$ preserves pullbacks (this
        can been seen by observing that the Grothendieck construction
        $\int_{\sA} \colon \pcat{\sA} \to \Cat/\sA$ has a left adjoint, and the
        forgetful functor $\Cat/\sA \to \Cat$ preserves pullbacks). In
        particular, it follows from \cref{paragr:internalgrpd} that we have an
        induced adjunction
        \[
          \Grpd(i_{\sA}) \colon \Grpd(\psh{\sA}) \rightleftarrows \Grpd(\Cat)
          \colon \Grpd(i_{\sA}^*).
        \]
        Note also that we have a canonical isomorphism $\Grpd(\psh{\sA})\simeq
        \pgrpd{\sA}$.
      \end{paragr}
      \begin{lemma}\label{lemma:createwe}
        Let $\sA$ be a small category. The functor \[\Grpd(i_{\sA}) \colon
          \pgrpd{\sA} \simeq \Grpd(\psh{\sA}) \to \Grpd(\Cat)\] preserves and
        reflects weak equivalences.
      \end{lemma}
      \begin{proof}
        Consider the commutative square (up to isomorphism) from
        \cref{lemma:comtensprod}
        \[
          \begin{tikzcd}[column sep=huge]
            \pgrpd{\sA} \simeq \Grpd(\psh{\sA})\ar[d,"N_*"']
            \ar[r,"\Grpd(i_{\sA})"] & \Grpd(\Cat) \ar[d,"N_*"]
            \\
            \psh{\Delta\times \sA} \simeq \Homu(\Delta^{\op},\psh{\sA})
            \ar[r,"{\Homu(\Delta^{\op},i_{\sA})}"'] & \pcat{\Delta}.
          \end{tikzcd}
        \]
        Let us say that a morphism $f$ of $\pcat{\Delta}$ is a weak equivalence
        if $\int_{\Delta}f$ is in $\cW$. Then, by definition, the vertical arrow
        on the right in the above square preserves and reflects weak
        equivalences, and it follows from \cref{lemma:intmult} that the bottom
        horizontal arrow of the above square also preserves and reflects weak
        equivalences. By \cref{cor:wenerve}, the vertical arrow on the left of
        the previous square preserves and reflects weak equivalences, and the
        desired conclusion follows at once.
      \end{proof}
      \begin{prop}\label{prop:grpdcatmodel}
        Let $\sA$ be a weak test category. The following adjunction
        \[
          \begin{tikzcd}[column sep=large]
            \Grpd(i_{\sA}) \colon \pgrpd{\sA} \simeq \Grpd(\psh{\sA})
            \ar[r,shift left] & \ar[l,shift left] \Grpd(\Cat) \colon
            \Grpd(i_{\sA}^*)
          \end{tikzcd}
        \]
        is a homotopical equivalence (\cref{def:homotopicalequiv}). In
        particular, it induces an adjoint equivalence after localization
        \[
          \begin{tikzcd}[column sep=large]
            \Ho(\pgrpd{\sA}) \ar[r,shift left,"\simeq"] & \ar[l,shift
            left,"\simeq"] \Ho(\Grpd(\Cat)).
          \end{tikzcd}
        \]
      \end{prop}
      \begin{proof}
        From \cref{lemma:createwe}, we know that $\Grpd(i_{\sA})$ preserves and
        reflects weak equivalences. Hence, by \cref{lemma:crithomotopicalequiv},
        we need to show that for every object $X=(X_0,X_1)$ of $\Grpd(\Cat)$,
        the co-unit
        \[
          \epsilon_{X} \colon \Grpd(i_{\sA}\circ i_{\sA}^*)(X) \to X
        \]
        is a weak equivalence of $\Grpd(\Cat)$. By definition, this means that
        we have to prove that $\int_{\Delta}N_{\ast}(\epsilon_{X})$ is in $\cW$.
        It follows from the fact that $i_{\sA}$ preserves pullbacks that
        $N_{\ast}(\Grpd(i_{\sA}\circ~i_{\sA}^*)(X))$ can be identified with the
        simplicial object in $\Cat$
        \[
          \begin{aligned}
            \Delta^{\op} &\to \Cat \\
            [n] &\mapsto
            i_{\sA}^*i_{\sA}(\underbrace{X_{1}\times_{X_0}\cdots\times_{X_0}X_1}_{n\text{
                times}})
          \end{aligned}
        \]
        and that $N_{\ast}(\epsilon_{X})$ can be identified with the morphism of
        simplicial objects of $\Cat$ which is the co-unit of the adjunction
        $i_{\sA}\dashv i_{\sA}^*$ level-wise
        \[
          i_{\sA}i_{\sA^*}(X_{1}\times_{X_0}\cdots\times_{X_0}X_1)\longrightarrow
          X_{1}\times_{X_0}\cdots\times_{X_0}X_1
        \]
        (cf.\ \cref{paragr:internalgrpd}). By hypothesis, $\sA$ is a weak test
        category and thus these morphisms are all in $\cW$. The conclusion
        follows then from \cref{ThomColimit}.
      \end{proof}
      \begin{rem}\label{rem:grpdcatmodel}
        In particular, when $\cW=\cW_{\infty}$, the previous proposition shows
        that $\Grpd(\Cat)$ models homotopy types. This is slightly surprising
        considering that groups internal to categories only model pointed
        connected homotopy $2$\nbd{}types
        \cite{loday1982spaces},\cite{maclane19503}. Somehow, restricting to
        groups internal to categories instead of groupoids internal to
        categories does not only amount to considering pointed connected
        objects. This is in contrast with what happens for simplicial sets:
        groupoids internal to simplicial sets models homotopy types and groups
        internal to simplicial sets models pointed connected homotopy types
        \cite{kan1958homotopy}.

        A hint of explanation comes from the fact that the functor $i_{\sA}
        \colon \psh{\sA} \to \Cat$ does not preserve products (even if it
        preserves pullbacks, as we have already seen), and in particular does
        not preserve groups objects. Hence, the strategy to prove
        \cref{prop:grpdcatmodel} cannot be adapted for group objects instead of
        groupoid objects.
      \end{rem}
      \begin{rem}
        The proof of \cref{prop:grpdcatmodel} does not depend on the theory of
        groupoidal test categories (as long as we define weak equivalences in
        $\Grpd(\psh{\sA})$ via the nerve as in \cref{sec:wenerve}). In
        particular, if we already know that $\Grpd(\psh{\Delta})$ models
        homotopy types (for example, by appealing to \cite[Theorem
        8.3]{crans1995quillen} or \cite[Theorem 10]{joyal1996homotopy}), then we
        obtain a second proof of the fact that $\Grpd(\psh{\sA})$ models
        homotopy type for any (weak) test category. Indeed, first we apply
        \cref{prop:grpdcatmodel} for $\sA=\Delta$ to deduce that $\Grpd(\Cat)$
        models homotopy types and then we apply the same result again for an
        arbitrary (weak) test category $\sA$. Note, however, that we do
        \emph{not} recover this way the other direction in the equivalence (b)
        of \cref{firstcompthm}.
      \end{rem}
      \bibliographystyle{plain} \bibliography{mybib}
    \end{document}